\documentclass[12pt]{article} \usepackage{amsmath}
\usepackage{amssymb} \usepackage[mathscr]{eucal}
\newcommand{\eqdef}{:=} 
\newcommand{\re}{\mathop{\tt Re}}
\newcommand{\finishproof}{\hfill $\Box$ \vspace{5mm}}

\newcommand{\p}{\partial} 

\newcommand{\At}{{\mathtt A}}

\newcommand{\D}{{\mathcal D}}
\newcommand{\Dt}{{\mathtt D}}
\newcommand{\Fk}{{\mathcal F}_{k}}
\newcommand{\Dl}{{\mathcal D}^{\ell}}
\newcommand{\hl}{H^{\ell}}
\newcommand{\C}{{\mathbb C}}
\newcommand{\K}{{\mathbb K}}
\newcommand{\T}{\mathbb{T}}
\newcommand{\N}{\mathbb{N}}
\newcommand{\Z}{\mathbb{Z}}
\newcommand{\R}{\mathbb{R}}

\newtheorem{Th}{Theorem}[section]
\newtheorem{Prop}[Th]{Proposition}
\newtheorem{Lemma}[Th]{Lemma}
\newtheorem{Coro}[Th]{Corollary}
\newtheorem{Def}[Th]{Definition}
\newtheorem{Rem}[Th]{Remark}

\numberwithin{equation}{section}

\begin{document}

\title{Analyticity of Riemannian exponential maps on ${\rm Diff}(\T)$}

\author{T. Kappeler\thanks{Supported in part by the Swiss National
Science Foundation, the programme SPECT, and the European
Community through the FP6 Marie Curie RTN ENIGMA
(MRTN-CT-2004-5652).}, E. Loubet \thanks{Supported by the Swiss
National Science Foundation }\,, P. Topalov}

\maketitle

\begin{abstract}
We study the exponential maps induced by Sobolev type right-invariant 
(weak) Riemannian metrics of order $k\ge1$ on the Lie group of smooth, 
orientation preserving diffeomorphisms of the circle. 
We prove that each of them defines an {\em analytic} Fr\'echet chart of the identity.
\end{abstract}

\section{Introduction}
The aim of this paper is to contribute towards a development of
Riemannian geometry for infinite dimensional Lie groups which has
attracted a lot of attention since Arnold's seminal paper \cite{Arn}
on hydrodynamics -- see \cite{AK}, \cite{CK2}, \cite{EM}, \cite{KM},
\cite{Kour}, \cite{Mis1}.  As a case study we consider the Lie group
${\mathcal D}\equiv{\rm Diff}({\mathbb T})$ of orientation preserving
$C^\infty $-diffeomorphisms of the circle ${\mathbb T} = {\mathbb
  R}/{\mathbb Z}$. According to Milnor \cite[$\S 9$]{Mil}, the group
$\D$, endowed with the $C^{\infty}$-Fr\'echet differential structure,
is \emph{not} analytic i.e., the mapping $\D\times\D\to\D$,
$(\varphi,\psi)\mapsto\varphi\circ\psi^{-1}$ is not analytic.
Nevertheless, it turns out that for a family of right-invariant weak
Riemannian metrics on $\D$, the corresponding Riemannian exponential
map defines an \emph{analytic} chart of the identity of $\D$ -- see
Theorem~\ref{Theorem 1.2} below.

The Lie group ${\mathcal D}$ and its Lie algebra come up in
hydrodynamics, playing the role of a configuration space for Burgers
equation and the Camassa-Holm equation \cite{Kour}, \cite{Mis1} (see
also \cite{KM}). The latter equation is a model for (one dimensional)
wave propagation in shallow water (cf.  \cite{FF}, and for a
derivation based on physical grounds \cite{CH}) that has several
interesting features which have been intensively studied in recent
years.

\medskip 

For any given integer $k \geq 0$, consider the scalar product $\langle
\cdot , \cdot \rangle _k : C^\infty ({\mathbb T}) \times C^\infty
({\mathbb T}) \rightarrow {\mathbb R}$,
\[ 
\langle u, v \rangle _k := \sum ^k_{j = 0} \int _{\mathbb T}
\partial ^j_x u \partial ^j_x v dx.
\]
It induces a $C^{\omega}_{F}$ (i.e., Fr\'echet analytic) weak, right-invariant, Riemannian
metric $\nu ^{(k)}$ on ${\mathcal D}$:
\[
\nu^{(k)}_{\varphi}(\xi,\eta):=
\langle(d_{\mathrm{id}}R_{\varphi})^{-1}\xi,(d_{\mathrm{id}}R_{\varphi})^{-1}\eta\rangle_{k},\quad\forall
\varphi\in{\mathcal D},\quad\text{and}\quad \forall\xi,\eta\in
T_{\varphi}\D
\]
where $R_{\varphi}:\D\rightarrow\D$ denotes the right translation
$\psi\mapsto\psi\circ\varphi$. The subscript $F$ in $C^{\omega}_{F}$
refers to the calculus in Fr\'echet spaces -- see
Appendix~\ref{Appendix:Frechet_spaces} where we collect some
definitions and notions of the calculus in Fr\'echet spaces. The
metric $\nu^{(k)}$ being weak means that the topology induced
by $\nu^{(k)}$ on the tangent space $T_{\varphi}\D$ at an arbitrary
point $\varphi$ in $\D$, is weaker than the Fr\'echet topology on
$C^{\infty}(\T)$.

\begin{Def}\hspace{-2mm}{\bf .}\label{definition 1.0}
For any given $T>0$, a $C^{2}_{F}$-smooth curve
$\varphi:[0,T]\rightarrow\D$, is called a geodesic with respect to
$\nu^{(k)}$, or $\nu^{(k)}$-geodesic for short, if it is a critical
point of the action functional within the class of $C^{2}_{F}$-smooth
variations $\gamma$ of $\varphi$ constrained to keep the end points fixed i.e., for
any $C^{2}_{F}$-smooth function
\[
\gamma:(-\varepsilon,\varepsilon)\times[0,T]\rightarrow\D,\,(s,t)\mapsto\gamma(s,t)
\]
such that $\gamma(0,t)=\varphi(t)$ and, $\gamma(s,0)=\varphi(0)$ and
$\gamma(s,T)=\varphi(T)$ for any $-\varepsilon<s<\varepsilon$, one has
\begin{subequations}
\begin{equation}\label{E_k'}
\frac{d}{ds} \Big\arrowvert_{s=0}\mathcal{E}_{k}^{T}(\gamma(s,\cdot))=0.
\end{equation}
where $\mathcal{E}_{k}^{T}$ denotes the action functional
\begin{equation}\label{E_k}
\mathcal{E}_{k}^{T}(\gamma(s,\cdot)):=\frac{1}{2}\int_{0}^{T}\nu_{\gamma(s,t)}^{(k)}(\dot\gamma(s,t),\dot\gamma(s,t))dt\,,
\end{equation}
\end{subequations}
and $\dot\gamma(s,t)=\p\gamma(s,t)/\p t$.
\end{Def}
The Euler-Lagrange equations \eqref{E_k'} on $\D$ for critical points
of the action functional $\mathcal{E}_{k}^{T}$ defined in \eqref{E_k}
are given by
\begin{equation}\label{ds'}
\left\{\begin{aligned}
{\dot\varphi}&=v\\
{\dot v}&=F_k(\varphi,v)
\end{aligned}\right.
\end{equation}
where
\begin{subequations}\label{FAB}
\begin{equation}
\label{F_k} F_k(\varphi , v):= R_{\varphi}\circ A^{-1}_k \circ B_k(v \circ\varphi ^{-1})
\end{equation}
in which
\begin{align}\label{A_k} 
A_k:= &\sum ^k_{j = 0} (-1)^j \partial ^{2j}_x 
\end{align}
and, for any smooth function $u$ in $C^{\infty}(\T)$,
\begin{eqnarray}\label{B_k} 
B_k(u)& := &- 2u' A_k u + A_k(uu') - uA_k u'\nonumber\\
&=& - 2u' A_k u + Q_k(u)
\end{eqnarray}
\end{subequations}
where $Q_k$ is a polynomial in $2k+1$ variables, and $Q_{k}(u)$ is
short for the function $Q_{k}(u,\p_{x}u,\ldots,\p_{x}^{2k}u)$. In the
sequel, we will refer to $Q_{k}(u)$ as the polynomial in $u,\partial_{x}u,\cdots,\partial^{2k}_{x}u$.
Here $(\cdot)^{\cdot}$ stands for $d/dt$ and $(\cdot)'$ for
$\p/\p_{x}$. Note that $t\mapsto\varphi(t)$ evolves in $\D$ whereas
$t\mapsto\dot\varphi(t)=v(t)$ is a vector field along $\varphi$ i.e.,
a section of $\varphi(\cdot)^{*}T\D$. In particular, $v(t)$ belongs to
$T_{\varphi(t)}\D$ for any $0\le t\le T$ for some $T>0$. We want to study the
following initial value problem
\begin{equation}
\label{1.3bis}
\left\{\begin{aligned}
    (\dot \varphi, \dot v) &= \left( v, F_k(\varphi , v) \right)\\
    (\varphi(0),v(0))&=(\mathrm{id},v_{0})\,.
\end{aligned}\right.
\end{equation}
It is easy to check that
\eqref{1.3bis} is equivalent to
\begin{align}
\label{1.5}
&\left\{\begin{aligned}
\dot\varphi&=u\circ\varphi\\
\varphi(0)&=\mathrm{id}
\end{aligned}\right.\\
\intertext{and}
\label{e:CH_k}
&\left\{\begin{aligned}
&A_k \dot u + u A_k u' + 2u' A_k u = 0\\
&u(0)=v_{0}
\end{aligned}\right.
\end{align}
where $t\mapsto u(t)=(d_{\mathrm{id}} R_{\varphi (t)})^{-1}\dot
\varphi (t)$ belongs to $T_{\mathrm{id}}{\mathcal D}$. The initial value
problems \eqref{1.3bis} and \eqref{e:CH_k} are, via \eqref{1.5}, two alternative descriptions
of the geodesic flow. The first corresponds to the Lagrangian
description i.e., tracking the moving point in ${\mathcal D}$ along
the section $\varphi(\cdot)^{*}T\D$ while the latter describes the
events in $T_{\mathrm{id}}{\mathcal D}$ from the Eulerian point of view
of a fixed observer.

\medskip

The case $k = 1$ is particularly interesting since the geodesic flow
\eqref{1.3bis} with $v=u\circ \varphi$ is a
re-expression of the Camassa-Holm equation \eqref{e:CH_k} -- see e.g.
\cite{KM},\cite{Kour}, and \cite{Mis2}. Indeed, the dynamical system
\eqref{1.3bis} can be used to study the
initial value problem for \eqref{e:CH_k} -- see e.g. \cite{DKT},
\cite{Mis2}.

\medskip

\begin{Th}\hspace{-2mm}{\bf .}\label{Theorem 1.1} 
  Let the integer $k\ge1$. Then, for any of the right-invariant metrics $\nu
  ^{(k)}$, there exists an open neighborhood $V^{(k)}$ of $0$ in $C
  ^\infty ({\mathbb T})$ such that, for any $v_0$ in $V^{(k)}$, there is
  a unique $\nu ^{(k)}$-geodesic, $(-2, 2) \rightarrow {\mathcal D}, t
  \mapsto \varphi (t;v_{0})$, issuing from the identity in the direction
  $v_0$. Moreover, the map
\[ 
(-2, 2)\times V^{(k)}\rightarrow{\mathcal D},\ (t,v_0)\mapsto \varphi (t; v_0)
\]
is Fr\'echet analytic.
\end{Th}

Theorem~\ref{Theorem 1.1} allows to define, for any given $k \geq 1$,
the Riemannian exponential map
\[ 
\mathrm{Exp}_k\big\arrowvert_{V^{(k)}} : V^{(k)} \rightarrow
{\mathcal D}, \quad v_0 \mapsto \varphi (1; v_0) .
\]

\begin{Th}\hspace{-2mm}{\bf .}\label{Theorem 1.2} 
For any integer $k\ge1$, there exists an open neighborhood
$\tilde{V}^{(k)}\subseteq V^{(k)}$ of $0$ in $C^\infty ({\mathbb T})$ and an open neighborhood $U^{(k)}$ of $\mbox{\rm id}$ in
${\mathcal D}$ so that
\[
\mathrm{Exp}_k\big\arrowvert_{\tilde{V}^{(k)}} : \tilde{V}^{(k)}\rightarrow U^{(k)},\quad
v_0\mapsto\varphi(1;v_{0})
\]
is a Fr\'echet bianalytic diffeomorphism.
\end{Th}

\medskip

{\it Related work:} Theorems~\ref{Theorem 1.1} -~\ref{Theorem 1.2}
improve the results of Constantin and Kolev in
\cite{CK1,CK2} where it was shown that the exponential map is a $C^1_F$-diffeomorphism. 

\medskip

{\it Method:} The approach used to establish Theorems~\ref{Theorem 1.1}
-~\ref{Theorem 1.2} is new. It consists in showing that the
analyticity of the geodesic flow as described in
Theorem~\ref{Th:local_existence}, Theorem~\ref{Th:extension},
Proposition~\ref{Proposition 6.1}, and Proposition~\ref{Proposition
  7.1}, can be obtained by an interplay of the structure of the
equations describing the geodesic flow from a Lagrangian perspective
and the structure of the Euler equations. In \cite{KLT}, we applied
our method to define and study exponential maps on the Lie group of
orientation preserving diffeomorphisms on the two dimensional torus
$\R^{2}/\Z^{2}$. Moreover, this new approach improves on the results in
\cite{CKKT} for the Virasoro group.

\medskip

{\it Organization of the paper:} Sections~ \ref{2. Group of
  diffeomorphism}, \ref{Sec:The vector field}, and \ref{Sec:Complex
  extension} are preliminary. In particular, we show in detail that the
vector fields in the equations defining the geodesic flows are
analytic in suitable spaces. Theorem~\ref{Theorem 1.1} is proved in
section~\ref{Sec:Exponential map}. To show Theorem~\ref{Theorem 1.2},
we use Theorem~\ref{Th:IFT} which is a version of the inverse function
theorem in a set-up for analytic maps between Fr\'echet spaces
discussed in Appendix~\ref{Appendix:Frechet_spaces}. In
section~\ref{8. Proof of Theorem}, we verify that the assumptions of
Theorem~\ref{Th:IFT} are satisfied in our situation.

\medskip

{\it Case $k=0$:} Our new technique does not apply in the case $k=0$;
the equation \eqref{ds'} with $k=0$, is no longer a dynamical system
(hence, we can not rely on the existence theorem of ODE's to establish
existence of geodesics in Sobolev spaces). Indeed, the crucial
difference with respect to the case $k\ge1$, lies in the fact that the
inverse of the operator $A_{k}$ defined in \eqref{A_k} is the identity
operator in the case $k=0$, hence it is \emph{not} regularizing. More
specifically, the expression on the r.h.s.  of the system \eqref{ds'}
is the map ${\mathcal D}^\ell \times H^\ell \rightarrow H^\ell \times
H^{\ell-1}$, $(\varphi,v)\mapsto(v,-2vv'/\varphi')$ which is
\emph{not} a vector field on $\Dl\times\hl$. Therefore the proof of
Theorem~\ref{Theorem 1.1} in the case $k=0$ has to be dealt with
differently. One way is to notice that the geodesic flow \eqref{ds'}
with $k=0$ is equivalent to the inviscid Burgers equation $\dot
u+3uu'=0$ (cf. \eqref{e:CH_k}) which can be solved implicitly (locally
in time) by the method of characteristics. However, in the case $k=0$,
Theorem~\ref{Theorem 1.2} does \emph{not} hold. An explicit counter
example is given in \cite{CK2}. For this reason, in the rest of the
paper, we will only concentrate on the case $k\ge1$.

\medskip

{\it Notation:} The notation we use is standard. In particular,
$H^s=H^s ({\mathbb T})=H^{s}(\T,\R)$ denotes the space of real valued
functions on ${\mathbb T}$ of Sobolev class $H^{s}$, and for $s \geq
1, W^{s , \infty }$ denotes the Banach space of continuous functions
$f : {\mathbb T} \rightarrow {\mathbb R}$ for which $\partial ^j_x f$
are in $L^{\infty}=L^{\infty}(\T)$ for every $0 \leq j \leq s $.
Further, for $s \geq 2$, ${\mathcal D}^s={\mathcal D}^s({\mathbb T})$
denotes the set of orientation preserving $C^{1}$-diffeomorphisms
$\varphi : {\mathbb T} \rightarrow {\mathbb T}$ of class $H^s $. It is
a Hilbert manifold modeled on $H^s$. The complexification of a real
vector space $X$ will be denoted by $X_\C$ i.e., $X_\C:= X\otimes\C$. 

\medskip

{\bf Acknowledgement:} It is a great pleasure for the first author to
thank P. Michor for valuable discussions.
\section{The group of diffeomorphisms of ${\mathbb T}$}
\label{2. Group of diffeomorphism}
For the convenience of the reader, we collect in this section the
results about the group of diffeomorphisms of $\T$ that we need. For a
general discussion on this topic see \cite{EM} and the references
therein.  Throughout this section, $s\ge2$. Denote by ${\mathcal D}^s$ 
the set of all orientation preserving
$C^{1}$-diffeomorphisms $\varphi : {\mathbb T} \rightarrow {\mathbb T}$ such that 
$\varphi '$ belongs to $H^{s - 1}$ i.e.,
\[
\mathcal{D}^{s}:=\{\varphi:\mathbb{T}\rightarrow\mathbb{T},\,\,\,
C^{1}-\mathrm{diffeomorphism}|\quad\varphi'>0,\quad\mathrm{and}\quad\varphi'\in
H^{s-1}\}.
\]
${\mathcal D}^s $ is in a natural way a Hilbert
manifold modeled by the Hilbert space $H^s $.  An atlas of
$\D^{s}$ can be described in terms of the lifts of $\varphi$ in
${\mathcal D}^s $. A lift of $\varphi$ is of the form
\[ 
{\mathbb R} \rightarrow {\mathbb R}, \quad x \mapsto x + f(x) 
\]
where $f$ is in $H^s$. The following two Hilbert charts form an atlas
of $\D^s$
\begin{subequations}
\begin{equation}\label{U_1}
{\mathfrak U}_1^s:= \{\varphi=\mathrm{id}+f|\;\; f \in
H^s\quad\mbox{and}\quad |f(0)| < 1 / 2 ; \ f' > - 1 \}
\end{equation}
and
\begin{equation}\label{U_2}
{\mathfrak U}_2^s:= \{\varphi=\mathrm{id}+f|\;\; f \in
H^s \quad\mbox{and}\quad\ 0
< f(0) < 1 ;\ f' > - 1 \} .
\end{equation}
\end{subequations}
(By a slight abuse of notation we denoted a lift of a
diffeomorphism $\varphi $ again by $\varphi $.)

\begin{Lemma}\hspace{-2mm}{\bf .}
\label{Lemma 2.1} Let $s\ge 2$. Then, for any $\varphi$ in ${\mathcal D} ^s $,
the inverse $\varphi ^{-1}$ of $\varphi $ is again in ${\mathcal D}^s $.
\end{Lemma}

{\it Proof.} Clearly, for any $\varphi$ in ${\mathcal D}^s $,
$\varphi ^{-1}$ is a $C^{1}$-diffeomorphism. Using that $(\varphi
^{-1})' = (1 / \varphi')\circ\varphi^{-1}$ and the fact that $H^r$ is
an algebra for any $r \geq 1$, one sees that $(\varphi ^{-1})'$ is in
$H^{s - 1}$, and hence that $\varphi ^{-1}$ belongs to ${\mathcal D}^s $.
\finishproof

\begin{Lemma}\hspace{-2mm}{\bf .}
\label{Lemma 2.2} Let $s \geq 2$. Then the following statements hold:
\begin{itemize}
\item[(i)] For any $u$ and $\varphi$ in $H^s$ and
$\mathcal{D}^s$ respectively, $u \circ \varphi$ is in $H^s$.
\item[(ii)] For any $\varphi$ in ${\mathcal D}^s $, the
  right-translation $R_\varphi : H^s \rightarrow H^s , \ u
  \mapsto u \circ \varphi $ is uniformly continuous on subsets $\mathcal{W}
  \subseteq {\mathcal D}^s $ satisfying
\begin{equation}
\label{2.0} 
\sup _{\varphi \in \mathcal{W}} (\| \varphi ^{-1} \| _{W^{1,\infty }} +
                \| \varphi \| _{H^s } ) < +\infty  .
\end{equation}
\item[(iii)] For any $u$ in $H^s$, the left-translation $L_u :
  {\mathcal D}^s \rightarrow H^s , \ \varphi \mapsto u \circ\varphi $
  is continuous.
\end{itemize}   
\end{Lemma}

{\it Proof.} $(i)$ First let us prove that for $u$ in $L^{2}$ and $\varphi$ in ${\mathcal
D}^2$, the composition $u \circ \varphi$ is in $L^{2}$. To see this, note that 
   \begin{equation}
     \label{2.1} \|u\circ\varphi\|_{L^{2}}^{2}=\int _{\mathbb T} |u(\varphi (x))|^2 dx = \int _{\mathbb T}
     |u(x)|^2 (\varphi ^{-1})' (x) dx \leq (\inf_{\T} \varphi ')^{-1} \| u\| ^2_{L^{2}}.
   \end{equation}

   Moreover, as $s\ge2$, by the Sobolev embedding theorem, $\varphi'$
   is bounded for any $\varphi$ in $\D^s$. Hence, the argument above
   shows that, for any $s\ge2$, given any $u$ in $H^{s}$ and any
   $\varphi$ in $\D^{s}$, $u\circ\varphi$ and $(u \circ \varphi )' =
   (u' \circ \varphi)\varphi'$ are in $L^{2}$. For any $2 \leq n \leq
   s $, we get by the chain and product rules
   \begin{equation}
\label{2.2manuscript}
 \partial ^n_x(u \circ \varphi ) = (u'\circ\varphi) \partial
      ^n_x \varphi + \sum ^n_{j = 2}  p_{n,j}(\varphi) [ (\partial ^j_x u)
      \circ \varphi]
   \end{equation}
where $p_{n,j}(\varphi)$ is a polynomial in $\p_x\varphi , \ldots ,
\partial ^{n + 1 - j}_x \varphi $. Now, for $s\ge 2$, the fact that
$u'$ is in $H^{s-1}\subseteq H^{1}\subseteq C^{0}$ implies that
$u'\circ\varphi$ is bounded so that
$(u'\circ\varphi)\p_{x}^{n}\varphi$ is in $H^{s-n}\subseteq L^{2}$.
Similarly, for every $2\le j\le n$, $(\p_{x}^{j}u)\circ\varphi$ are
seen to be in $L^{2}$, and $p_{n,j}(\varphi)$ is bounded. But then, as
$\p_{x}^{n}(u\circ\varphi)$ is in $L^{2}$ for every $2\le n\le s$,
and since $u\circ\varphi$ and $(u\circ\varphi)'$ were shown to be in
$L^{2}$ too, we conclude that, for any $s\ge2$, $u$ in $H^{s}$ and $\varphi$ in $\D^{s}$
imply that $u\circ\varphi$ is in $H^{s}$ as asserted.

\medskip

$(ii)$ Let $(u_m)_{m \geq 1}$ be a sequence in $H^s $ which converges
to $u$ in $H^s $. By \eqref{2.1}, one has
\[
\| (u_m - u)\circ\varphi\|_{L^{2}}\leq\|\varphi
^{-1}\|^{1/2}_{W^{1,\infty }}\| u_m - u \|_{L^{2}},
\]
and
   \[ 
\| ((u_m - u) \circ \varphi )' \| _{L^{2}} \leq \| \varphi \| ^{1/2}
      _{W^{1,\infty }} \| u_m - u \| _{H^{1}} .
   \]
   To estimate the derivatives of order $2 \leq n \leq s $, we
   use formula \eqref{2.2manuscript}.

By the Sobolev embedding theorem, 
   \[ \| [(u_m - u)' \circ \varphi] \partial ^n_x \varphi \| _{L^{2}} \leq
      \| u_m - u \| _{W^{1, \infty }} \| \varphi \| _{H^n} \leq C \| u_m - u\|
      _{H^s } \| \varphi \| _{H^s },
   \]
   where $C>0$ denotes a constant, and for every $2 \leq j \leq n$,
   \[ \| p_{n,j}(\varphi) [\partial ^j_x (u_m - u) \circ \varphi] \|
   _{L^{2}} \leq [q_{n,j} (\| \varphi \| _{H^s })]^{1/2} \| \varphi
   ^{-1}\| ^{1/2}_{W^{1,\infty }} \| u_m - u \| _{H^s }
   \]
where $q_{n,j}$ is a polynomial in one variable with positive
coefficients. It then follows from \eqref{2.2manuscript} and the above
estimates that, for any $0\le j\le s$,
\[
\|\p_{x}^{j}[R_{\varphi}(u_{m}-u)]\|_{L^{2}}=O(\|u_{m}-u\|_{H^{s}})
\]
for $\varphi$ in sets $\mathcal{W}\subseteq\D^{s}$ satisfying
\eqref{2.0}. This establishes uniform continuity of the right
translation within the set $\mathcal{W}$.

\medskip

$(iii)$ Let $(\varphi _m)_{m \geq 1} \subseteq {\mathcal D}^s $ be
convergent to $\varphi$ in ${\mathcal D}^s $. First we will assume
that $u$ belongs to $H^{s + 2}$. As $s \geq 2$, $u$ is Lipschitz
continuous i.e., for any $x$ in $\T$,
   \[ |u \circ \varphi _m(x) - u \circ \varphi (x)| \leq \| u \| _{W^{1, \infty }}
      |\varphi _m(x) - \varphi (x)|.
   \]
This implies
   \[ \|u \circ \varphi _m - u \circ \varphi \| _{L^{2}} \leq \| u \| _{W^{1, \infty }}
      \|\varphi _m - \varphi \|_{L^{2}} .
   \]
   Similarly, as $(u\circ\varphi_{m})'=
   (u'\circ\varphi_{m})\varphi'_{m}$, one can estimate
   $\|(u\circ\varphi_{m})'-(u\circ\varphi)'\|_{L^{2}}$ by
\begin{align*}
& \|(u'\circ\varphi_m)(\varphi'_m-\varphi')\|_{L^{2}}+ 
\|(u'\circ\varphi_m-u'\circ\varphi)\varphi'\| _{L^{2}}\\
\leq&\|u\|_{W^{1,\infty}}\|(\varphi_m-\varphi)'\|_{L^{2}}+\|\varphi\|_{W^{1,\infty }} 
\|u\|_{W^{2,\infty }}\|\varphi_m-\varphi\|_{L^{2}}\\
\leq&C\,\mathrm{max}(1,\|\varphi\|_{W^{1,\infty}})\|u\|_{H^{3}}\|\varphi_{m}-\varphi\|_{L^{2}},
\end{align*}
where in the last step we used the estimate $\| u\| _{W^{1, \infty}} \leq\| u\| _{W^{2, \infty }} \leq C\| u\| _{H^3}$, 
which follows by the Sobolev embedding theorem.

\medskip

Next, as in the proof of $(ii)$, we use \eqref{2.2manuscript} to estimate the higher order derivatives 
($2\leq n\le s$). First,
\[
\begin{array}{l}
\|\partial^{n}_{x}(u\circ\varphi_m - u\circ\varphi)\|_{L^{2}}\le
\|(u'\circ\varphi_m)\partial^{n}_{x}\varphi_{m}-(u'\circ\varphi)\partial^{n}_{x}\varphi\|_{L^{2}}\\
\quad\quad\quad\quad\quad\quad\quad\quad\quad
+\sum^{n}_{j = 2}\|p_{n,j}(\varphi_{m})[(\partial^{j}_{x}u)\circ\varphi_{m}]-p_{n,j}(\varphi)
[(\partial^{j}_{x}u)\circ\varphi]\|_{L^{2}}.
\end{array}
\]
Based on previous estimates, the norm
$\|(u'\circ\varphi_{m})\partial^{n}_{x}\varphi_{m}-(u'\circ\varphi)\partial^{n}_{x}\varphi\|_{L^{2}}$ 
is bounded by 
\[
\|u'\circ\varphi_{m}-u'\circ\varphi\|_{L^{\infty}}\|\partial^{n}_{x}\varphi_{m}\|_{L^{2}}+\|u\|_{W^{1,\infty}}
\|\p_{x}^{n}\varphi_{m}-\p_{x}^{n}\varphi\|_{L^{2}},
\]
while $\|
p_{n,j}(\varphi_{m})[(\partial^{j}_{x}u)\circ\varphi_{m}]-p_{n,j}(\varphi)
[(\partial^{j}_{x}u)\circ\varphi]\|_{L^{2}}$ is bounded by
\[
\|(p_{n,j}(\varphi_{m})-p_{n,j}(\varphi))[(\partial^{j}_{x}u)\circ\varphi_{m}]\|_{L^{2}}+\|
p_{n,j}(\varphi)[(\partial^{j}_{x}
u)\circ\varphi_{m}-(\partial^{j}_{x} u)\circ\varphi]\|_{L^{2}}.
\]
One then shows that $\| \partial ^n_x(u \circ \varphi _m - u \circ\varphi ) \| _{L^{2}}\rightarrow0$ as
$m\rightarrow+\infty$ using the estimates
\[
\begin{array}{l} 
\|u'\circ\varphi_m-u'\circ\varphi\|_{L^\infty}\|\partial^n_x\varphi_m\|_{L^{2}}\leq 
\|u\|_{W^{2, \infty}}\sup_{m\ge 1}\{\|\varphi_m\|_{H^n}\}\|\varphi_m-\varphi\|_{H^{1}},\\
\end{array}
\]
\[
\begin{array}{l}
\|u'\circ\varphi\|_{L^\infty }\|\partial^n_x\varphi_m-\partial^n_x\varphi\|_{L^{2}}\leq 
\|u\|_{W^{1,\infty }}\|\varphi _m - \varphi \| _{H^n}, \\
\end{array}
\]
and
\[
\begin{array}{l}
\|(p_{n,j}(\varphi_{m})-p_{n,j}(\varphi))[(\partial^{j}_{x}u)\circ\varphi_{m}]\|_{L^{2}}\le\\
\le\| \partial ^j_x u\| _{L^{2}} \| \varphi ^{-1}_m \|^{1/2}_{W^{1,\infty }}
\| p_{n,j}(\varphi_{m}) - p_{n,j}(\varphi)\|_{L^\infty }\\
\le\| u\| _{H^s } \sup_{m\ge1}\{\| \varphi ^{-1}_m\|^{1/2}_{W^{1, \infty }}\}
r_{n,j}(\|\varphi _m-\varphi\|_{W^{s-1,\infty }})
\end{array}
\]
where $r_{n,j}$ is a polynomial in one variable with positive coefficients, and
\begin{align*} 
\| p_{n,j}(\varphi)[(\partial^{j}_{x}u)\circ\varphi_{m}-(\partial^{j}_{x} u)\circ\varphi]\|_{L^{2}}
\le\| p_{n,j}(\varphi) \| _{L ^\infty }\| \partial ^j_x u\|_{W^{1,\infty }}\| \varphi _m - \varphi \|_{L^{2}}\\
\le q_{n,j} (\| \varphi \| _{W^{n - 1, \infty }})\| u\|_{W^{j+1,\infty}} \| \varphi _m - \varphi \|_{L^{2}}\\
\le C q_{n,j} (\| \varphi \| _{W^{s - 1, \infty }})\| u\|_{H^{s + 2}} \| \varphi _m - \varphi \| _{L^{2}}
\end{align*}
   where the last inequality results from the fact that, for every
   $2\le j\le s$, $\|\p_{x}^{j} u\|_{W^{1,\infty}}\le\|
   u\|_{W^{j+1,\infty}}\le\|u\|_{W^{s+1,\infty}}\le C\| u\|_{H^{s+2}}$
   by the Sobolev embedding theorem, and where $q_{n,j}$ is a
   polynomial in one variable with positive coefficients.  Altogether,
   since $u$ was assumed to be in $H^{s+2}$, we have shown that
   $(L_{u}-L_{u_{m}})\circ\varphi\rightarrow0$ as $m\rightarrow
   +\infty$ in $H^{s}$.

\medskip

The general case, where $u$ is merely in $H^{s}$, is now obtained by a
limiting argument. Indeed, as $H^{s+2}$ is densely embedded in
$H^{s}$, we may choose a sequence $(u_i)_{i \geq 1} \subseteq
H^{s + 2}$ so that $u_i \rightarrow u$ in $H^s $ as
$i\rightarrow +\infty$. Then
   \[ \| u \circ \varphi _m - u \circ \varphi \| _{H^s } \leq \| u \circ
      \varphi _m - u_i \circ \varphi _m\| _{H^s } + \| u_i \circ \varphi _m
      - u_i \circ \varphi \| _{H^s } + \| u_i \circ \varphi - u \circ \varphi
      \| _{H^s }.
   \]
   Note that
   $\mathrm{sup}_{m\ge1}(\|\varphi_{m}^{-1}\|_{W^{1,\infty}}+\|\varphi_{m}^{-1}\|_{H^{s}})<+\infty$
   as $\varphi_{m}\rightarrow\varphi$ in $\D^{s}$ for $s\ge2$. Hence,
   we have by $(ii)$ that, for any given $\varepsilon>0$, there exists
   $i_{0}=i_{0}(\varepsilon)\geq 1$ so that, for any $i \geq i_0$, and
   any $m\geq 1$, the first and third factors on the r.h.s. of the
   last display are no larger than $\varepsilon/3$. Finally, by the
   argument above, there exists $m_{0}(\varepsilon)\ge1$ such that for
   every $m\ge m_{0}$ the middle term is smaller than $\varepsilon/3$,
   and we are done.  \finishproof

\begin{Prop}\hspace{-2mm}{\bf .}\label{Proposition 2.3} 
  For any $s\ge2$ and any $r$ in $\Z_{\ge0}$, the composition $H^s \times
  {\mathcal D}^{s+r}\rightarrow H^{s}$, $(u,\varphi ) \mapsto u \circ \varphi $ is
  $C^{r}$-smooth.
\end{Prop}
{\it Proof.} Consider first the case $r=0$. It is to prove that, for any sequences $(u_m)_{m \geq 1} \subseteq
H^s $ and $(\varphi _m)_{m \geq 1} \subseteq {\mathcal D}^s $
such that $u_m \rightarrow u$ in $H^{s}$ and $\varphi_m
\rightarrow\varphi$ in $\D^{s}$ as $m\rightarrow +\infty$,
the sequence $(u_m \circ \varphi _m)_{m \geq 1}$ converges in $H^s $ to $u\circ \varphi $. Indeed,
\[ 
\|u_m \circ \varphi _m - u \circ \varphi \| _{H^s } \leq \| (u_m - u)\circ \varphi _m \| _{H^s } + 
\| u \circ \varphi _m - u \circ\varphi \| _{H^s } .
\]
As $\varphi _m$ converges to $\varphi$ in ${\mathcal D}^s $ and
$s \geq 2$, one has that $ \sup _m \| \varphi ^{-1}_m \| _{W^{1,\infty }} < +\infty$. 
Hence one can apply Lemma~\ref{Lemma 2.2} $(ii)$ to conclude that $\| (u_m - u) \circ \varphi _m \| _{H^s}\rightarrow 0$. 
By Lemma~\ref{Lemma 2.2} $(iii)$, $\| u \circ\varphi _m - u \circ \varphi \| _{H^s }\rightarrow 0$ as
$m\rightarrow +\infty$. The proof of the case $r\ge1$ is similar and
is left to the reader.
\finishproof

\begin{Prop}\hspace{-2mm}{\bf .}\label{Proposition 2.4} 
  For every $s\ge2$ and any $r$ in $\Z_{\ge0}$, the map ${\mathcal D}^{s+r}
  \rightarrow {\mathcal D}^s$, $\varphi\mapsto\varphi^{-1}$ is
  $C^{r}$-smooth.
\end{Prop}

{\it Proof.} As above, we will give the proof of the case $r=0$, and
leave the case $r\ge1$ to the reader. Let $(\varphi _m)_{m \geq 1}$ be
a sequence in ${\mathcal D}^s $ with $\varphi _m \rightarrow
\varphi $ in ${\mathcal D}^s $. By Lemma~\ref{Lemma 2.1},
$\varphi_{m}^{-1}$ and $\varphi^{-1}$ are in $\D^{s}$. It is to prove
that $\varphi ^{-1}_m \rightarrow \varphi ^{-1}$ in ${\mathcal D}^s
$.  To this end, write
\[ 
\varphi ^{-1}_m - \varphi ^{-1}= (\mathrm{id} - \varphi ^{-1}\circ\varphi _m) \circ \varphi ^{-1}_m .
\]
By Lemma~\ref{Lemma 2.2} $(iii)$, $\varphi ^{-1} \circ\varphi_{m}\rightarrow\varphi ^{-1} \circ \varphi = \mathrm{id}$ in
$\D^{s}$ as $m\rightarrow +\infty$. As $\sup _m (\| \varphi ^{-1}_m \|_{H^s } + \| \varphi _m \| _{W^{1, \infty }}) < +\infty $, it
then follows from Lemma~\ref{Lemma 2.2} $(ii)$ that $(\varphi^{-1}\circ\varphi_{m}-\mathrm{id})\circ\varphi_{m}^{-1}\rightarrow0$
as $m\rightarrow +\infty$.
\finishproof
   
\begin{Rem}
From Proposition~\ref{Proposition 2.3} and
Proposition~\ref{Proposition 2.4} it follows that the composition and the inverse maps,
$\D\times\D\rightarrow\D$ respectively $\D\rightarrow\D$, are
$C^{\infty}_{F}$-smooth, hence $\D$ is a Lie group. Its Lie algebra
$T_{\mathrm{id}}\D$ can be canonically identified with $C^{\infty}(\T,\R)$, with bracket given by $[u,v]=uv'-u'v$.
\end{Rem}

\section{The vector field $\Fk$}\label{Sec:The vector field}
In the present section, let the integer $k \geq 1$, and $\ell \geq \ell _k := 2k + 2$.
For any $(\varphi , v)$ in ${\mathcal D}^\ell \times H^\ell$, consider
\begin{equation}\label{ds} 
{\mathcal F}_k(\varphi , v):= (v, F_k(\varphi , v)),
\end{equation}
where $F_{k}$ is defined in \eqref{FAB}. 
It follows from Lemma~\ref{Lemma 2.1} and Lemma~\ref{Lemma 2.2} that,
for any $(\varphi,v)$ in $\D^\ell\times H^\ell$, the r.h.s. of 
\eqref{ds} is well-defined and belongs to the space $H^\ell \times H^\ell$.
In particular, \eqref{ds} defines a {\em dynamical system} (ODE) on $\D^\ell\times H^\ell$.
Introducing 
\[
{\cal A}_k : \D^\ell\times H^\ell \to\D^\ell\times H^{\ell-2k} ,\;\;(\varphi,v)\mapsto(\varphi,R_\varphi\circ A_k\circ R_{\varphi^{-1}}v)
\]
and
\[
{\cal B}_k : \D^\ell\times H^\ell \to\D^\ell\times H^{\ell-2k} ,\;\;(\varphi,v)\mapsto(\varphi,R_\varphi\circ B_k\circ R_{\varphi^{-1}}v)
\]
with $A_k$ and $B_k$ as in \eqref{A_k}, respectively \eqref{B_k}, we
can write
\[
F_k={\mathcal Proj}_2\circ{\cal A}_k^{-1}\circ{\cal B}_k
\]
where ${\mathcal Proj}_2$ is the projection onto the second component
$(\varphi,v)\mapsto v$, and ${\cal A}_k^{-1}$ is the inverse of $
{\cal A}_k$ described in the following proposition.
\begin{Prop}\hspace{-2mm}{\bf .}\label{Lem:crucial}
Let $k\ge 1$, and $\ell\ge\ell_k=2k+2$. Then 
\begin{itemize}
\item[(i)] the map
\begin{subequations}
\label{defcalAkinv}
\begin{equation}\label{A_k-smooth}
{\cal A}_k : \D^\ell\times H^\ell \to\D^\ell\times H^{\ell-2k} ,\;\;(\varphi,v)\mapsto(\varphi,R_\varphi\circ A_k\circ R_{\varphi^{-1}}v)
\end{equation}
is a bianalytic diffeomorphism with inverse given by the map
\begin{equation}\label{A_k-inverse}
{\cal A}_k^{-1} : \D^\ell\times H^{\ell-2k} \to\D^\ell\times H^{\ell} ,\;\;(\varphi,v)\mapsto(\varphi,R_\varphi\circ A_k^{-1}\circ R_{\varphi^{-1}}v)\,.
\end{equation}
\end{subequations}
\item[(ii)] The map
\begin{equation}\label{B_k-smooth}
{\cal B}_k : \D^\ell\times H^\ell\to\D^\ell\times H^{\ell-2k} ,\;\;(\varphi,v)\mapsto(\varphi,R_\varphi\circ B_k\circ R_{\varphi^{-1}}v)
\end{equation}
is analytic.
\end{itemize}
As a consequence,
\begin{itemize}
\item[(iii)]the vector field
  \begin{equation}\label{e:decomposition} {\mathcal F}_k
    :\D^\ell\times H^\ell\to H^\ell\times H^\ell ,\;\;
    (\varphi,v)\mapsto {\cal A}_k^{-1}\circ{\cal B}_k (\varphi,v)
\end{equation}
is analytic on the Hilbert manifold $\D^\ell\times H^\ell$.
\end{itemize}
\end{Prop}

To prove Proposition~\ref{Lem:crucial} we need two auxiliary lemmas.

\begin{Lemma}\hspace{-2mm}{\bf .}\label{Lem:auxiliary}
Let $k\ge 1$, and $\ell\ge\ell_k=2k+2$. Then, for any $\varphi$ in
$\D^\ell$ and $v$ in $H^\ell$, and, for any $1\le n\le 2k$, the following statements hold:
\begin{itemize}
\item[(i)]$\p_{x}(\varphi^{-1})\circ\varphi=1/\varphi'$, and for
  $n\ge2$, $(\partial_{x}^{n}\varphi^{-1})\circ\varphi$ is a
  polynomial in $1/\varphi'$,
  $\partial_{x}\varphi,\ldots,\partial_{x}^{n}\varphi$ with integer
  coefficients.
\item[(ii)]
  $\mathtt{D}^{n}(\varphi,v):=R_{\varphi}\circ\p_{x}^{n}\circ
  R_{\varphi^{-1}}v=(\p_x^n(v\circ\varphi^{-1}))\circ\varphi=\sum_{j=1}^{n}p_{n,j}(\varphi)\p_x^jv$
  where $p_{n,j}(\varphi)$ is a polynomial in $1/\varphi'$,
  $\p_x\varphi,\ldots,\p_x^{n+1-j}\varphi$ with integer coefficients.
  In particular, $p_{n,n}(\varphi)=1/(\varphi')^{n}$.
\end{itemize}
\end{Lemma}
{\em Proof.}
The proof follows by a straightforward application of the chain rule
(cf. discussion following \eqref{2.2manuscript}).\finishproof

Recall that $H^{s}_{\C}:=H^{s}(\T,\C)$ is the complexification of
$H^{s}$.
\begin{Lemma}\hspace{-2mm}{\bf .}\label{Lem3.3}
  For $s\ge 1$, let $W^{s}_{\C}$ denote the open subset $W^{s}_\C:=\{
  f \in H^{s}_{\C} \; : \; f(x)\ne 0\;\forall x\in\T\} $. Then, the
  map $W^{s}_{\C}\to H^{s}_{\C}$, $f\mapsto 1/f$ is analytic.
\end{Lemma}
{\em Proof.}
Let $f$ in $W^{s}_{\C}$, and $U_{\epsilon,\C}(f)$ be the neighborhood
\[
U_{\epsilon,\C}(f)=\{f-g\;|\;g\in H^{s}_{\C},\;
;\;\|g\|_{H^s_\C}<\epsilon\}
\]
with $\epsilon>0$ so small that $\|g/f\|_{H^{s}_{\C}}<1$. Such a
choice is possible since, $H^s_\C$ being a Banach algebra for $s\ge1$,
$\|g/f\|_{H^{s}_{\C}}\le C\|g\|_{H^{s}_{\C}} \|1/f\|_{H^{s}_{\C}}$ so
that it suffices to pick $0<\epsilon<1/(C \|1/f\|_{H^{s}_{\C}})$.
Then, $1/(f-g)$ can be written in terms of a series
\[
\frac{1}{f-g}=\frac{1}{f}\left(1+\frac{g}{f}+\left(\frac{g}{f}\right)^2+\ldots\right)
\]
which converges uniformly in $U_{\epsilon,\C}(f)$ to an element in $H^s_\C$.\finishproof

\begin{Coro}\hspace{-2mm}{\bf .}\label{Corollary3.4}
  For any $s\ge 2$, the map $\D^s\to H^{s-1}$, $\varphi\mapsto
  1/\varphi'$ is analytic.
\end{Coro}
{\em Proof.}
The map $\D^{s}\rightarrow H^{s-1}$, $\varphi\mapsto1/\varphi'$ is
the composition of the linear map $\D^{s}\rightarrow
H^{s-1}$, $\varphi\mapsto\varphi'$, and the analytic map
$H^{s-1}\rightarrow H^{s-1}$, $\varphi'\mapsto1/\varphi'$.\finishproof

\medskip

{\em Proof of Proposition~\ref{Lem:crucial}.} $(i)$ By direct
computation, one sees that the map defined in \eqref{A_k-inverse} is
indeed the inverse of \eqref{A_k-smooth}. In particular, this shows
that $\mathcal{A}_{k}$ is bijective. By Lemma~\ref{Lem:auxiliary}, the
definition \eqref{A_k} of $A_k$, and that of $\Dt^{n}(\varphi,v)$ (cf.
Lemma~\ref{Lem:auxiliary} $(ii)$), we have that
\begin{align}
R_\varphi\circ A_k\circ
R_{\varphi^{-1}}v&=v+\sum_{j=1}^{k}(-1)^{j}\Dt^{2j}(\varphi,v)\nonumber\\
&=v+\sum_{j=1}^{2k}q_{2k,j}(\varphi)\p_x^jv\label{defAtk1}
\end{align}
where $q_{2k,j}(\varphi)$ is a polynomial in $1/\varphi'$,
$\p_x\varphi,\ldots,\p_x^{2k+1-j}\varphi$ with integer
coefficients. Note that in view of Lemma~\ref{Lem:auxiliary} $(ii)$, $q_{2k,2k}(\varphi)=(-1)^{k}/(\varphi')^{2k}$.

By Corollary~\ref{Corollary3.4}, the map
\[
\D^\ell\times H^\ell \to (H^{\ell-2k}
)^{4k+2}\,,\quad(\varphi,v)\mapsto(1/\varphi',\p_x\varphi,...,\p_x^{2k}\varphi,v,\p_x
v,...,\p_x^{2k}v)
\]
is analytic and, since for $\ell\ge\ell_{k}=2k+2$, $H^{\ell-2k}$ is a
Banach algebra, we conclude that the r.h.s. of
\eqref{defAtk1} is analytic and hence that $\mathcal{A}_{k}$ is
analytic.  Moreover, for any $(\varphi_0,v_0)$ in $\D^\ell\times
H^\ell$, the differential $d_{(\varphi_0,v_0)}{\cal A}_k : H^\ell
\times H^\ell\to H^\ell \times H^{\ell-2k} $ is of the form
\begin{equation}\label{e:differential}
  d_{(\varphi_0,v_0)}{\cal A}_k(\delta\varphi,\delta v)=
  \left(
\begin{array}{cc}
\delta\varphi&0\\
\Lambda(\delta\varphi)&R_{\varphi_0}\circ A_k\circ R_{\varphi_0^{-1}}\delta v
\end{array}
\right)
\end{equation}
where 
\[
\Lambda : H^\ell\to H^{\ell-2k}\,,\quad\mbox{and}\quad
R_{\varphi_0}\circ A_k\circ R_{\varphi_0^{-1}} : H^\ell \to
H^{\ell-2k}
\]
are bounded linear maps.  As the latter map is invertible, the open
mapping theorem implies that it is a linear isomorphism.  Hence,
$d_{(\varphi_0,v_0)}{\cal A}_k$ is a linear isomorphism and, by the
inverse function theorem, the map defined in \eqref{A_k-smooth} is a
local $C^\omega$-diffeomorphism and since we have seen that
$\mathcal{A}_{k}$ is bijective, assertion $(i)$ follows.  The proof of
item $(ii)$ is similar to the proof of the analyticity of ${\cal A}_k$
in part $(i)$.\finishproof

\section{The complex analytic extension of ${\mathcal F}_k$}\label{Sec:Complex extension}
As in the previous section, let $k\ge 1$, and $\ell\ge\ell_k=2k+2$.
Denote by ${\mathfrak U}_{1,\C}^\ell$ the complexification of the
Hilbert chart ${\mathfrak U}_1^\ell$ defined in \eqref{U_1},

\begin{equation}\label{U_1-complex} {\mathfrak U}_{1,\C}^\ell:=
  \{\varphi=\mathrm{id}+f\,|\;\;f\in H^\ell_\C \, ;\,|f(0)| < 1 /
  2;\,\re\,[(\varphi')^{2k-3}]>0\}\,.
\end{equation}

(The condition $\re\,[(\varphi')^{2k-3}]>0$ will be used in the proof
of Proposition~\ref{Proposition 6.1}.)  It follows from
Lemma~\ref{Lem:auxiliary} and Lemma~\ref{Lem3.3} that, for any $1\le
n\le 2k$, the map
\begin{subequations}
\label{defDtnRC}
\begin{equation}\label{e:D^s-real}
{\mathfrak U}_{1}^\ell\times H^{\ell}\to
H^{\ell-n},\quad (\varphi,v)\mapsto
\Dt^{n}(\varphi,v):=R_{\varphi}\circ\p_{x}^{n}\circ R_{\varphi^{-1}}v
\end{equation}

can be extended to an analytic map

\begin{equation}\label{e:D^s-complex}
{\mathfrak U}_{1,\C}^\ell\times H^\ell_\C\to H^{\ell-n}_\C\,,\quad (\varphi,v)\mapsto
\Dt^{n}_{\C}(\varphi,v).
\end{equation}
\end{subequations}

As a consequence, the map (cf. \eqref{defAtk1})

\begin{subequations}
\label{defAtkAtkC}
\begin{equation}
\label{defAtk}
{\mathfrak U}_{1}^\ell\times H^{\ell}\to
H^{\ell-2k}\,,\quad(\varphi,v)\mapsto\At_{k}(\varphi,v):=v+\sum_{j=1}^{k}(-1)^{j}\Dt^{2j}(\varphi,v)
\end{equation}

has an analytic extension 
\[ {\mathfrak U}_{1,\C}^\ell\times H^{\ell}_\C \to
H^{\ell-2k}_\C\,,\quad(\varphi,v)\mapsto\At_{k,\C}(\varphi,v):=v+\sum_{j=1}^{k}(-1)^{j}\Dt_{\C}^{2j}(\varphi,v)\,.
\]
Note that the latter is of the form
\begin{equation}
\label{defAtkC}
\At_{k,\C}(\varphi,v)=v+\sum_{j=1}^{2k}q_{2k,j}(\varphi)\p_{x}^{j}v
\end{equation}
\end{subequations}
where $q_{2k,j}(\varphi)$ is a polynomial in
$1/\varphi',\p_{x}\varphi,\ldots,\p_{x}^{2k+1-j}\varphi$. Further, we
introduce the analytic map (the complexification of \eqref{A_k-smooth})

\begin{equation}\label{e:A_kC^-1}
{\mathcal A}_{k,\C} : {\mathfrak U}^\ell_{1,\C}\times H^\ell_\C\to {\mathfrak U}^\ell_{1,\C}\times H^{\ell-2k}_\C,\,
(\varphi,v)\mapsto\left(\varphi,\At_{k,\C}(\varphi,v)\right)\,.
\end{equation}

Analogously, Lemma~\ref{Lem:auxiliary} and Lemma~\ref{Lem3.3} imply
that the map \eqref{B_k-smooth} can be analytically extended to the
map

\begin{equation}\label{e:B_k-complex}
{\mathcal B_{k,\C}} : {\mathfrak U}_{1,\C}^\ell\times H^\ell_\C\to{\mathfrak U}_{1,\C}^\ell\times H^{\ell-2k}_\C\,.
\end{equation}

\begin{Lemma}\hspace{-2mm}{\bf .}\label{Lem:analytic_extension1}
For any $k\ge1$ and $\ell\ge\ell_k$, there exists a complex neighborhood $U_{\ell,k;\C}$ of $\mathrm{id}$ 
in ${\mathfrak U}_{1,\C}^\ell$ such that \eqref{A_k-inverse} can be extended to a bianalytic diffeomorphism

\begin{equation}\label{e:A_kC}
{\mathcal A}_{k,\C}^{-1} : U_{\ell,k;\C}\times H^{\ell-2k}_\C\to U_{\ell,k;\C}\times H^{\ell}_\C
\end{equation}

with inverse ${\mathcal A}_{k,\C}\big|_{U_{\ell,k;\C}\times
  H^{\ell}_\C}$.
\end{Lemma}

{\em Proof.}  By Proposition~\ref{Lem:crucial}, for any
$\ell\ge\ell_k$, there exists a complex neighborhood
$U_{\ell,k;\C}\times W_{\ell,\C}^{(k)}$ of $(\mathrm{id},0)$ in
${\mathfrak U}^{\ell}_{1,\C}\times H^\ell_\C$ such that ${\mathcal
  A}_{k,\C}\big|_{U_{\ell,k;\C}\times W_{\ell,\C}^{(k)}}$ is a
bianalytic diffeomorphism onto its image.  In particular, cf.
\eqref{e:A_kC^-1}, for any given $\varphi$ in $U_{\ell,k;\C}$, the map
$\At_{k,\C}(\varphi,\cdot\,) : W^{(k)}_{\ell,\C}\to H^{\ell-2k}_\C$
extends, by linearity, to a bounded map

\[
\At_{k,\C}(\varphi,\cdot\,) : H^{\ell}_\C\to H^{\ell-2k}_\C\,.
\]

As $\At_{k,\C}(\varphi,\cdot\,)$ is a bianalytic diffeomorphism near
the origin, it is in fact bijective. Hence, by the open mapping
theorem, $\At_{k,\C}(\varphi,\cdot\,)$ is a linear isomorphism for any
$\varphi$ in $U_{\ell,k;\C}$ and hence

\begin{equation}\label{e:loc3}
{\mathcal A}_{k,\C}\big|_{U_{\ell,k;\C}\times H^{\ell}_\C} : 
U_{\ell,k;\C}\times H^{\ell}_\C\to U_{\ell,k;\C}\times H^{\ell-2k}_\C
\end{equation}

is bijective. Moreover, from \eqref{e:A_kC^-1}, one sees that, for any
$(\varphi_{0},v_{0})$ in $U_{\ell,k;\C}\times H^{\ell}_\C$, the
differential $d_{(\varphi_0,v_0)}{\cal A}_{k,\C}$ of the map in
\eqref{e:loc3} is of the form
\eqref{e:differential}, and hence, a linear isomorphism. Altogether, it
follows that \eqref{e:loc3} is a bianalytic
diffeomorphism.\finishproof

\begin{Lemma}\hspace{-2mm}{\bf .}\label{Lem:analytic_extension}
For any $k\ge1$ and any $\ell\ge\ell_k$, the neighborhood $U_{\ell,k;\C}$
in Lemma~\ref{Lem:analytic_extension1} can be chosen to be of the form
$U^{(k)}_{\ell,\C}:=U_{\ell_{k},k;\C}\cap{\mathfrak U}_{1,\C}^\ell$.
\end{Lemma}

{\em Proof.} For any $\ell\ge\ell_k$, formula \eqref{e:A_kC^-1}
defines an analytic map from $U_{\ell,\C}^{(k)}\times H^\ell_\C$ to
$U_{\ell,\C}^{(k)}\times H^{\ell-2k}_\C$.  By
Lemma~\ref{Lem:analytic_extension1}, for $\ell=\ell_k$, this map is a
bianalytic diffeomorphism. Hence, it is injective for any
$\ell\ge\ell_k$. The proof of
Lemma~\ref{Lem:analytic_extension1} shows that it suffices to prove that, for any $\varphi$ in
$U_{\ell,\C}^{(k)}$, 
\[
{\mathcal A}_{k,\C}(\varphi,\cdot\,) :
H^\ell_\C\to H^{\ell-2k}_\C\,,\quad v\mapsto \At_{k,\C}(\varphi,v)\,,
\]
is onto. Indeed, for any $\varphi$ in $U_{\ell,\C}^{(k)}$ and any $h$
in $H^{\ell-2k}_\C$, the equation for $v$, $\At_{k,\C}(\varphi,v)=h$,
is by \eqref{defAtkC} the ODE
\[
v+\sum_{j=1}^{2k}q_{2k,j}(\varphi)\p_x^jv=h\,.
\]
Now, first observe that for $\varphi$ in $U_{\ell,\C}^{(k)}\subseteq
U_{\ell_{k},\C}^{(k)}$, the linear operator
$\At_{k,\C}(\varphi,\cdot\,)$ maps $H^{\ell_{k}}_{\C}$ to
$H^{\ell_{k}-2k}_{\C}$, and since the latter is a linear isomorphism,
it follows that (for any $\varphi$ and $h$ as above) this equation has
a unique solution $v$ in $H^{\ell_{k}}_{\C}$.  Finally, as
$q_{2k,j}(\varphi)$ is in $H^{\ell-2k+j-1}_\C$ and
$q_{2k,2k}(\varphi)=(-1)^{k}/(\varphi')^{2k}$ does nowhere vanish by \eqref{U_1-complex}, it then follows that
$v$ is in $H^\ell_\C$.  \finishproof

From now on, we choose $U_{\ell,k;\C}$ in
Lemma~\ref{Lem:analytic_extension1} to be given by
\[
U^{(k)}_{\ell,\C}:=U_{\ell_{k},k,\C}\cap{\mathfrak U}_{1,\C}^\ell\,.
\]
Recall that the maps ${\mathcal A}_{k,\C}$ and ${\mathcal B}_{k,\C}$ described in \eqref{e:A_kC^-1} respectively
\eqref{e:B_k-complex} are analytic on ${\mathfrak U}^\ell_{1,\C}\times H^\ell_\C$.
By Lemma~\ref{Lem:analytic_extension1}, ${\mathcal A}_{k,\C}^{-1}$ is analytic on
$U_{\ell,\C}^{(k)}\times H^{\ell-2k}_\C$ (cf.  Lemma~\ref{Lem:analytic_extension}). Hence, the vector field
\begin{equation}\label{e:F_k-analytic}
{\mathcal F}_{k,\C} : U_{\ell,\C}^{(k)}\times H^\ell_\C\to
H^\ell_\C\times H^\ell_\C\,,\quad (\varphi,v)\mapsto(v,F_{k,\C}(\varphi,v)):=\mathcal{A}_{k,\C}^{-1}\circ\mathcal{B}_{k,\C}(\varphi,v)\,.
\end{equation}
is analytic; in fact the analytic extension of the vector field
defined in \eqref{e:decomposition}. We will study the properties of
the dynamical system corresponding to ${\mathcal F}_{k,\C}$ on
$U_{\ell,\C}^{(k)}$ i.e.,
\begin{equation}\label{ds-complex}
\left\{\begin{aligned}
{\dot\varphi}&=v\\
{\dot v}&=F_{k,\C}(\varphi,v)\,.
\end{aligned}\right.
\end{equation}
As, for any $k\ge 1$, $(\mathrm{id},0)$ is a zero of ${\mathcal
  F}_{k,\C}$ (and hence an equilibrium solution of
\eqref{ds-complex}), one gets from \cite[Theorem 10.8.1]{Dieudonne}
and \cite[Theorem 10.8.2]{Dieudonne} the following result.

\begin{Th}\hspace{-2mm}{\bf .}\label{Th:local_existence} 
  Let $k\ge 1$, and $\ell\ge\ell_{k}=2k+2$. Then there exists an open
  neighborhood $V_{\ell,k;\C}$ of $0$ in $H^\ell_\C $ so that, for any
  $v_{0}$ in $V_{\ell,k;\C}$, the initial value problem for
  \eqref{ds-complex} with initial data
  $(\varphi(0),v(0))=(\mathrm{id},v_{0})$ has a unique analytic
  solution
\begin{subequations}
\label{curves}
\begin{equation}
\label{curvevarphiv}
(-2,2)\rightarrow U_{\ell,\C}^{(k)}\times H^\ell_\C,\;t\mapsto(\varphi(t;v_{0}),v(t;v_{0})) .
\end{equation}
Moreover, the flow map,
\begin{equation}
\label{curveini}
(-2,2)\times V_{\ell,k;\C}\rightarrow U_{\ell,\C}^{(k)}\times
H^\ell_\C,\; (t,v_0)\mapsto(\varphi(t; v_0),v(t; v_0))
\end{equation}
\end{subequations}
is analytic.   
\end{Th}
\begin{Rem}
\label{remc1}
In fact, \cite[Theorem 10.8.1]{Dieudonne} and \cite[Theorem
10.8.2]{Dieudonne} imply that \eqref{curveini} is a $C^{1}$-map over
$\C$, and thus that the map \eqref{curveini} is analytic.
\end{Rem}
\begin{Rem}
\label{Remark 4.2}\hspace{-2mm}{\bf .} Theorem~\ref{Th:local_existence} does not exclude that
$\bigcap _{\ell\geq\ell_k}V_{\ell,k;\C}=\{0\}$.  This possibility is
ruled out by Theorem~\ref{Th:extension} given in the next section.
\end{Rem}

\section{The exponential map and its analytic extension}\label{Sec:Exponential map}
As in the previous sections, let $k \geq 1$, and set $\ell_k:=2k+2$.
By \eqref{curveini} of Theorem~\ref{Th:local_existence}, the
Riemannian exponential map
\[
\mathrm{Exp}_{k,\ell_k}: V_{\ell_{k},k;\C}\cap H^{\ell_{k}}\rightarrow
\D^{\ell_{k}},\; v_0\mapsto\varphi(1; v_0)
\]
admits an analytic extension
\begin{equation}
\label{e:Exp}
\mathrm{Exp}^{\C}_{k,\ell_k}: V_{\ell_{k},k;\C}\rightarrow U_{\ell_{k},\C}^{(k)},\; v_0\mapsto\varphi(1; v_0)\,. 
\end{equation}
Set 
\[
V^{(k)}_{\ell_{k},\C}:=V_{\ell_{k},k;\C}\,.
\]
Noting that
$d_{0}\mathrm{Exp}^\C_{k,\ell_k}=\mathrm{Id}_{H^{\ell_{k}}_{\C}}$, it
then follows from the inverse function theorem that, by shrinking the
neighborhoods $V_{\ell_{k},\C}^{(k)}$ and $U_{\ell_{k},\C}^{(k)}$ if
necessary, one can ensure that the mapping \eqref{e:Exp} is a
bianalytic diffeomorphism. This will be tacitly assumed in the
remaining of the paper. In this section, we study the restriction of
$\mathrm{Exp}^\C_{k,\ell_k}$ to $V^{(k)}_{\ell_{k},\C}\cap
C^\infty(\T,\C)$.

\vspace{0.5cm}

\noindent{\em A priori relation:}
For a while let us study the equation \eqref{ds'} instead of its
analytic extension \eqref{ds-complex}. Consider the curve
$u=v\circ\varphi^{-1}$ where $t\mapsto(\varphi(t),v(t))$ is an
arbitrary $C^{1}$-solution of \eqref{ds'} in
$(U_{\ell,\C}^{(k)}\cap{\mathfrak U}_1^\ell)\times H^\ell$ on some
nontrivial time interval $(-T,T)$. By
Proposition~\ref{Proposition 2.3} and Proposition~\ref{Proposition
  2.4}, 
\[
(-T,T)\to H^\ell\subseteq H^{\ell-1}\,,\quad t\mapsto
u(t)=v(t)\circ\varphi^{-1}(t)\,,
\]
is a $C^1$-curve in $H^{\ell-1}$.  Moreover, it satisfies
\eqref{e:CH_k}. Our aim is to derive, for any $-T<t<T$, a formula for $(A_k
u(t))\circ\varphi(t)$ which will be used to study regularity
properties of the exponential map.  By the chain rule,
\[ 
[(A_ku)\circ\varphi ]^{\cdot}=(A_k{\dot u})\circ\varphi +
[(A_k u)'\circ\varphi]\dot\varphi .
\]
As $v=\dot\varphi=u\circ\varphi$, the initial value problem \eqref{e:CH_k} then leads to
\[
\left\{
\begin{aligned}
&[(A_k u)\circ\varphi ]^{\cdot}+2(u'\circ\varphi)[(A_k u)\circ \varphi]=0\\
&(A_k u(0))\circ\varphi(0)=A_{k}v_{0}.
\end{aligned}
\right.
\]
By definitions \eqref{e:D^s-real} and \eqref{defAtk},
$(A_{k}u)\circ\varphi=\At_{k}(\varphi,v)$, and
$u'\circ\varphi=\Dt^{1}(\varphi,v)$. Hence, the latter can be
rewritten as
\[
\left\{
\begin{aligned}
  &[\At_k(\varphi,v)]^{\cdot}+2\Dt^1(\varphi,v) [\At_k(\varphi,v)]=0\\
  &\At_k(\varphi(0),v(0))=A_{k}v_{0}.
\end{aligned}
\right.
\]
where, for $-T<t<T$, $t\mapsto \At_k(\varphi(t),v(t))$ is of class $H^{\ell-2k}$ and
$t\mapsto \Dt^1(\varphi(t),v(t))$ evolves in $H^{\ell-1}$. Both of
these curves are $C^1$-smooth.
Solving the latter equation one gets
\begin{equation}\label{e:loc1} 
\At_k(\varphi(t),v(t))=e^{-2\int^t_0 \Dt^1(\varphi(\tau),v(\tau))\,d\tau}A_k v_0.
\end{equation}
On the other hand, differentiating of
\eqref{1.5} with respect to $x$ yields
\[
(\varphi')^\cdot= \Dt^1(\varphi,v)\varphi'\,.
\]
Since $\varphi (0)=\mathrm{id}$, one obtains
\begin{equation} 
\varphi'(t)=e^{\int^t_0 \Dt^1(\varphi(\tau),v(\tau))\,d\tau}\,.\label{e:phi_x*}
\end{equation}
Hence, \eqref{e:loc1} can be rewritten as
\[
\At_k(\varphi(t),v(t))=A_kv_0/(\varphi'(t))^2\,.
\]
Let
\[
(U_{\ell,\C}^{(k)}\cap{\mathfrak U}_1^\ell)\times
H^\ell\rightarrow
H^{\ell-2k},\quad(\varphi,v)\mapsto I_{k}(\varphi,v):=\At_{k}(\varphi,v)(\varphi')^{2}
\]
then, the above identity shows that the function
\[
I_{k}(\varphi(t),v(t))=\At_{k}(\varphi(t),v(t))(\varphi'(t))^2
\]
is independent of $t$, and is equal to $A_{k}v_{0}$. As a consequence, the
derivative $\mathcal{L}_{{\mathcal F}_k}(I_{k})$ of $I_{k}$ in the
direction ${\mathcal F}_k$ vanishes on the open set
$(U_{\ell,\C}^{(k)}\cap{\mathfrak U}_1^\ell)\times H^\ell$.

\medskip

Now, let us return to the analytic extension \eqref{ds-complex} of
\eqref{ds'}.  First note that
\begin{equation}\label{e:the_integral} 
  U_{\ell,\C}^{(k)}\times H^\ell_\C\rightarrow
  H^{\ell-2k}_{\C}\,,\quad (\varphi,v)\mapsto I_{k,\C}(\varphi,v):=\At_{k,\C}(\varphi,v)(\varphi')^2\,,
\end{equation}
is an analytic extension of $I_k$. Further, ${\mathcal F}_{k,\C}$
analytically extends ${\mathcal F}_k$, and the derivative
$\mathcal{L}_{{\mathcal F}_{k,\C}}(I_{k,\C})$ of $I_{k,\C}$ in the
direction ${\mathcal F}_{k,\C}$, analytically extends
$\mathcal{L}_{{\mathcal F}_k}(I_{k})$ to $U_{\ell,\C}^{(k)}\times
H^\ell_\C$ so that
\[
\mathcal{L}_{{\mathcal
    F}_{k,\C}}(I_{k,\C})\big|_{(U_{\ell,\C}^{(k)}\cap{\mathfrak
    U}_1^\ell)\times H^\ell}=\mathcal{L}_{{\mathcal F}_k}(I_{k})\,.
\]
As $\mathcal{L}_{{\mathcal F}_k}(I_{k})$ vanishes on
$(U_{\ell,\C}^{(k)}\cap{\mathfrak U}_1^\ell)\times H^\ell$, it then
follows that $\mathcal{L}_{{\mathcal F}_{k,\C}}(I_{k,\C})=0$
everywhere in $U_{\ell,\C}^{(k)}\times H^\ell_\C$, -- see e.g.
\cite[Proposition 6.6]{BS2}. In other words, \eqref{e:the_integral} is
a conserved quantity for the solutions of \eqref{ds-complex} on
$U_{\ell,\C}^{(k)}\times H^\ell_\C$.

\begin{Lemma}\hspace{-2mm}{\bf .}\label{Lem:relation}
For any $(\varphi,v)$ in $U^{(k)}_{\ell,\C}\times H^{\ell}_{\C}$, and
any $1\le j\le 2k$, the following relation holds
\begin{equation}
\label{complexndiffv}
\partial_x^j v=\Dt^j_\C(\varphi,v)(\varphi')^j+\Dt^1_\C(\varphi,v)\partial_x^j\varphi+\ldots
\end{equation}
where $\Dt^j_\C(\varphi,v)$ is the analytic extension of \eqref{e:D^s-real} to $U_{\ell,\C}^{(k)}\times H^\ell_\C$, and 
$\ldots $ stand for a polynomial in the variables $1/\varphi'$, $\p_{x}\varphi$,..., $\partial_x^{j-1}\varphi$ and
$\p_{x}v$,...,$\partial_x^{j-1}v$. 
\end{Lemma}

{\em Proof.}
First let us consider $(\varphi,v)$ in $(U_{\ell,\C}^{(k)}\cap{\mathfrak
    U}_1^\ell)\times H^\ell$. Then $\varphi$ is in $\Dl$ and
  $u:=v\circ\varphi^{-1}$ is well-defined in $H^{\ell}$, and thus we can
  write $v=u\circ\varphi$. By differentiation we get
\begin{subequations}
\begin{equation}
\label{5.6bis}
v'=(u'\circ\varphi)\varphi'\,.
\end{equation}
As $u'\circ\varphi=R_{\varphi}\circ\p_{x}\circ R_{\varphi^{-1}}v$, one
has by definition \eqref{e:D^s-real} that
\begin{equation}
\label{v'}
v'=\Dt^{1}(\varphi,v)\varphi'\,.
\end{equation}
By the analyticity of $\Dt_{\C}^{1}$,
the identity \eqref{v'} continues to hold for $(\varphi,v)$ in
$U_{\ell,\C}^{(k)}\times H^{\ell}_{\C}$ i.e.,
\[
v'=\Dt^{1}_{\C}(\varphi,v)\varphi'\,.
\]
For $j\ge2$, we argue similarly i.e., given $(\varphi,v)$ in
$(U_{\ell,\C}^{(k)}\cap{\mathfrak U}_1^\ell)\times H^\ell$, we differentiate \eqref{5.6bis} $(j-1)$ times to get
\begin{equation}\label{5.7bis}
\partial_x^j v=\Dt^j(\varphi,v)(\varphi')^j+\Dt^1(\varphi,v)\partial_x^j\varphi+\ldots
\end{equation}
\end{subequations}
where $\ldots$ stand for a polynomial in $1/\varphi'$,
$\p_{x}\varphi$,..., $\partial_x^{j-1}\varphi$ and
$\p_{x}v$,...,$\partial_x^{j-1}v$. Finally, \eqref{5.7bis} extends by analyticity
to $U_{\ell,\C}^{(k)}\times H^{\ell}_{\C}$ leading to
\eqref{complexndiffv}.\finishproof

Now, assume that $t\mapsto(\varphi(t),v(t))$ is a solution of
\eqref{ds-complex} in $C^{1}((-2,2),U_{\ell,\C}^{(k)}\times
H^\ell_\C)$. Then, by Lemma~\ref{Lem:relation}, the curve
$t\mapsto (A_k-1)\varphi(t)=\sum^k_{j = 1}(-1)^j\partial^{2j}_x\varphi(t)$
satisfies the inhomogeneous transport equation
\begin{equation}\label{e:inhomogeneous_equation} 
\left\{
\begin{aligned}
&((A_k-1)\varphi)^\cdot-\Dt^1_\C(\varphi,v) (A_k-1)\varphi=\At_{k,\C}(\varphi,v)(\varphi')^{2k}+g_{2k-1}(\varphi,v)\\
&(A_k-1) \varphi(0)=0
\end{aligned}
\right.
\end{equation}
where $g_{2k-1}(\varphi,v)$ is a polynomial (with constant
coefficients) in $1/\varphi'$, $\p_{x}\varphi$,...,
$\partial_x^{2k-1}\varphi$ and $v$, $\p_{x}v$,...,
$\partial_x^{2k-1}v$. (For convenience, we consider $(A_{k}-1)\varphi$
instead of $A_{k}\varphi$ so that the initial value problem
\eqref{e:inhomogeneous_equation} involves periodic functions only.)
Integrating \eqref{e:inhomogeneous_equation} by the method of
variation of parameters, and using \eqref{e:phi_x*} to write the final
expression in compact form, we get that, for any $-2<t<2$,
\[
  (A_{k}-1)\varphi(t)=\varphi'(t)\int^{t}_{0}
    \left[\At_{k,\C}(\varphi(\tau),v(\tau))(\varphi'(\tau))^{2k-1}+\dfrac{g_{2k-1}(\varphi(\tau),v(\tau))}{\varphi'(\tau)}\right]d\tau.
\]
By \eqref{e:the_integral} which when evaluated at the solution of
\eqref{ds-complex} is equal to $A_{k}v_{0}$, we then have that, for any
$-2<t<2$,
\begin{equation}\label{e:crucial} 
(A_k-1)\varphi(t)-\varphi'(t)\left(\int^t_0(\varphi'(\tau))^{2k-3}d\tau\right)
A_k v_0=\varphi'(t)\int^{t}_{0}\rho_{2k-1}(\varphi(\tau),v(\tau))\,d\tau 
\end{equation}
where 
\begin{equation}
\label{rho}
U_{\ell,\C}^{(k)}\times H^\ell_\C\to
H^{\ell-2k+1}_\C\,,\quad(\varphi,v)\mapsto\rho_{2k-1}(\varphi,v):=g_{2k-1}(\varphi,v)/\varphi'
\end{equation}
is analytic by Lemma~\ref{Lem3.3}. As we will explain in
detail later on, the a priori relation \eqref{e:crucial} plays a
fundamental role in the proofs of Theorem~\ref{Th:extension} and
Theorem~\ref{Theorem 1.2} below.

\begin{Th}\hspace{-2mm}{\bf .}\label{Th:extension}
  Let $k\ge1$, and $\ell\geq\ell_k=2k+2$. Then, for any $v_0$ in
  $V^{(k)}_{\ell,\C}:=V^{(k)}_{\ell_{k},\C}\cap H^\ell_\C$, there
  exists a unique solution of \eqref{ds-complex} in
  $C^1((-2,2),U_{\ell,\C}^{(k)}\times H^\ell_\C)$ with initial data
  $(\mathrm{id},v_0)$. Moreover, the flow map
\[
(-2,2)\times V_{\ell,\C}^{(k)}\rightarrow U_{\ell,\C}^{(k)}\times
H^\ell_\C,\; (t,v_0)\mapsto(\varphi(t; v_0),v(t; v_0))
\]
is analytic.
\end{Th} {\em Proof.}
We argue by induction with respect to $\ell\ge\ell_{k}$. For
$\ell=\ell_k$ the statement follows from Theorem
\ref{Th:local_existence} since, by definition,
$V^{(k)}_{\ell_{k},\C}=V_{\ell_{k},k;\C}$. Assume that the statement
is true for any given $\ell>\ell_{k}$ i.e., that given any $v_{0}$ in
$V^{(k)}_{\ell,\C}$, there exists a unique solution of
\eqref{ds-complex} in $C^1((-2,2),U_{\ell,\C}^{(k)}\times H^\ell_\C)$
with initial data $(\mathrm{id},v_{0})$ and, in addition, that the
flow map $(-2,2)\times V_{\ell,\C}^{(k)}\rightarrow
U_{\ell,\C}^{(k)}\times H^\ell_\C$ is analytic. Then, by \eqref{rho},
the r.h.s. of \eqref{e:crucial} is in $H^{\ell-2k+1}_\C$.  Now, let
$v_0$ be in $V^{(k)}_{\ell+1,\C}$ and, for $-2<t<2$, let
\[
t\mapsto\zeta(t):=(\varphi(t),v(t))\in U_{\ell,\C}^{(k)}\times
H^\ell_\C
\]
be the corresponding solution of \eqref{ds-complex} issuing from
$(\mathrm{id},v_{0})$. In particular, $t\mapsto\zeta(t)$ satisfies the
integral equation
\begin{equation}\label{e:integral_equation}
\zeta(t)=(\mathrm{id},v_{0})+\int_0^t{\mathcal F}_{k,\C}(\varphi(\tau),v(\tau))\,d\tau\,.
\end{equation}
As $v_0$ belongs to $V^{(k)}_{\ell+1,\C}$ and $k\ge1$, $A_{k}v_{0}$
and hence the second term on the l.h.s. of \eqref{e:crucial} are in
$H^{\ell-2k+1}_\C$. Altogether, it then follows from \eqref{e:crucial}
that $t\mapsto (A_k-1)\varphi(t)$ is a $C^1$-curve evolving in
$H^{\ell-2k+1}_\C$ for $-2<t<2$.  Hence, as $v=\dot\varphi$,
$t\mapsto\zeta(t)=(\varphi(t),v(t))$ is a continuous curve in
$U_{\ell+1,\C}^{(k)}\times H^{\ell+1}_\C$. By the analyticity of the
map \eqref{e:F_k-analytic},
it then follows that the integrand in \eqref{e:integral_equation} is a
continuous function of $\tau$ with values in $H^{\ell+1}_\C\times
H^{\ell+1}_\C$. Finally, the integral equation
\eqref{e:integral_equation} implies that $\zeta$ is a solution of
\eqref{ds-complex} in $C^1((-2,2),U_{\ell+1,\C}^{(k)}\times
H^{\ell+1}_\C)$ with initial data $(\mathrm{id},v_0)$. The second
statement of the theorem follows by combining \cite[Theorem
10.8.1]{Dieudonne} and \cite[Theorem 10.8.2]{Dieudonne} (cf. Remark~\ref{remc1}).\finishproof

\medskip

{\em Proof of Theorem~\ref{Theorem 1.1}.} Theorem~\ref{Theorem 1.1} is
an immediate consequence of Theorem~\ref{Th:extension}. Indeed, for
any $k\ge1$, 
\begin{equation}
\label{5.10bis}
V^{(k)}:=V_{\C}^{(k)}\cap C^{\infty}(\T,\R)
\end{equation}
satifies the properties stated in Theorem~\ref{Theorem
  1.1}.\finishproof

\medskip

Theorem~\ref{Th:extension} allows to define the exponential map.
Recall that, for any $\ell\geq\ell_k$,
\begin{subequations}
\begin{equation}
\label{defVkellC}
V^{(k)}_{\ell,\C}:=V_{\ell_{k},\C}^{(k)}\cap H^\ell_\C\,.
\end{equation}
where $V_{\ell_{k},\C}^{(k)}=V_{\ell_{k},k;\C}$. (Note that for
$\ell\geq\ell_k+1$, $V^{(k)}_{\ell,\C} $ might not coincide with the
neighborhood $V_{\ell,k;\C}$ introduced in
Theorem~\ref{Th:local_existence}.) Introduce
\begin{align}
V^{(k)}_\C&:=V_{\ell_{k},\C}^{(k)}\cap C^\infty({\mathbb T},\C)=
\displaystyle\bigcap_{\ell\ge\ell_{k}}V_{\ell,\C}^{(k)}\label{5.11bis}\\
U^{(k)}_\C&:=U_{\ell_{k},\C}^{(k)}\cap C^\infty({\mathbb T},\C)=
\displaystyle\bigcap_{\ell\ge\ell_{k}}U_{\ell,\C}^{(k)}\,.\nonumber
\end{align}
\end{subequations}
By Theorem~\ref{Th:extension}, for any $\ell\ge\ell_{k}$, the
restriction $\mathrm{Exp}^\C_{k;\ell}$ of $\mathrm{Exp}^\C_{k,\ell_k}$
to $V^{(k)}_{\ell,\C} $ takes values in $U_{\ell,\C}^{(k)}$. Moreover,
$\mathrm{Exp}^\C_{k,\ell} : V^{(k)}_{\ell,\C}\rightarrow
U_{\ell,\C}^{(k)}$ is analytic. Hence, the restriction
$\mathrm{Exp}^\C_k$ of $\mathrm{Exp}^\C_{k,\ell_k}$ to $V^{(k)}_{\C}$
takes values in $U^{(k)}_{\ell_{k},\C}\cap C^{\infty}(\T,\C)$ i.e.,
\[ 
\mathrm{Exp}^\C_k : V^{(k)}_{\C}\rightarrow U^{(k)}_\C\,,\quad
v_{0}\mapsto\varphi(1;v_{0})\,.
\]
\section{A Fr\'echet analytic chart of $\mathrm{id}$ in $\D(\T)$}\label{8. Proof of Theorem} 
Theorem~\ref{Theorem 1.2} states that the exponential map
$\mathrm{Exp}_k\big\arrowvert_{V^{(k)}} : V^{(k)} \rightarrow
{\mathcal D}$ that is defined using Theorem~\ref{Theorem 1.1},
can be used to define an analytic chart of the identity in $\D$.

\medskip

{\it Proof of Theorem~\ref{Theorem 1.2}.} By definitions
\eqref{5.10bis} and \eqref{5.11bis}, $V^{(k)}=V_{\C}^{(k)}\cap
C^{\infty}(\T,\R)$, and the map $\mathrm{Exp}^\C_k$, defined at the
end of the previous section, is the analytic extension of
$\mathrm{Exp}_{k}$. We want to apply to $\mathrm{Exp}^{\C}_{k}$ the
inverse function theorem in Fr\'echet spaces, Theorem~\ref{Th:IFT}.
Fix $k \geq 1$. To match the notation of this theorem, we write, for
any integer $n\ge0$, $\ell:=\ell_{k}+n$ ($\ell_{k}=2k+2$), and define
\[
X_n:=H^\ell_\C;\quad Y_n:=H^\ell_\C,\quad\quad\mbox{and} \quad\quad
V_n:= V^{(k)}_{\ell,\C};\quad
U_n:=\mathrm{Exp}^\C_{k,\ell}(V^{(k)}_{\ell,\C} )
\]  
where $V^{(k)}_{\ell,\C} $ is defined in \eqref{defVkellC}, and
$\mathrm{Exp}^\C_{k, \ell }$ is the exponential map introduced in
section~\ref{Sec:Exponential map}.  
Further, let $f:= \mathrm{Exp}^{\C}_{k,\ell _k} : V_{0} \rightarrow U_{0}$. 
It follows from our construction that $f$ is a $C^{1}$-diffeomorphism
and hence item $(a)$ of Theorem~\ref{Th:IFT} is verified.
Assumption (b) holds in view of Theorem~\ref{Th:extension}, whereas
items $(c)$ and $(d)$ hold, respectively, by
Proposition~\ref{Proposition 6.1}, and Proposition~\ref{Proposition 7.1} below. 
Hence, Theorem~\ref{Theorem 1.2} follows from Theorem~\ref{Th:IFT}.  
\finishproof

\medskip
  
It remains to show the two propositions used in the proof above.

\begin{Prop}\hspace{-2mm}{\bf .}\label{Proposition 6.1} 
Let $k \geq 1$, and $\ell_k=2k+2$. Then, for any $\ell=\ell_k+n$,
$n\ge0$, and any $v_0$ in $V^{(k)}_{\ell_{k},\C}$,
\[ 
\mathrm{Exp}^{\C}_{k,\ell_{k}}(v_{0})\in
U_{\ell,\C}^{(k)}\;\;\;\Longrightarrow \;\;\;v_{0}\in
V^{(k)}_{\ell,\C} .
\]
\end{Prop}
{\it Proof.} 
Let $(\varphi(\cdot;v_0),v(\cdot;v_0))$ denote the solution of
\eqref{ds-complex} with initial data $(\mathrm{id},v_0)$ with $v_{0}$
in $V^{(k)}_{\ell_{k},\C}$, and suppose that $(A_{k}-1)\varphi(1;v_{0})$
belongs to $H^{\ell-2k}_\C$. We will show that $v_{0}$ has to be in
$H^\ell_\C$.  For $\ell=\ell_{k}$, the result holds by construction.
Now, inductively, assume that $v_{0}$ is in $V^{(k)}_{\ell,\C}$ for
any given $\ell>\ell_{k}$. By Theorem~\ref{Th:extension}, this
solution actually lies in $C^1((-2,2),U_{\ell,\C}^{(k)}\times
H^\ell_\C)$. This implies that the r.h.s. of \eqref{e:crucial}, when
evaluated at $t=1$, is in $H^{\ell-2k+1}_\C$. Moreover, as
$U_{\ell,\C}^{(k)}\subseteq{\mathfrak U}_{1,\C}^\ell$, we have by
definition \eqref{U_1-complex} that the factor
$\varphi'(1)\int_{0}^{1}(\varphi'(\tau))^{2k-3}\,d\tau$ does not
vanish and is in $H^{\ell-1}_\C\subseteq H^{\ell-2k+1}_\C$ since
$k\ge1$. Altogether, it follows from \eqref{e:crucial} that
$A_{k}v_{0}$ lies in $H^{\ell-2k+1}_\C$, and hence that $v_{0}$ is in
$H^{\ell+1}_\C$.  \finishproof

By Theorem \ref{Th:extension}, $(-2,2)\times V_{\ell,\C}^{(k)}\rightarrow U_{\ell,\C}^{(k)}\times
H^\ell_\C,\; (t,v_0)\mapsto(\varphi(t; v_0),v(t; v_0))$ is analytic. Then, the variation $\delta v_0$ in
$H^\ell_\C$ of the initial data $v_0$ in $V_{\ell,\C}^{(k)}$ induces the variation of
$t\mapsto(\varphi(t;v_{0}),v(t;v_{0}))$
\[
t\mapsto(\delta\varphi(t),\delta
v(t)):=\dfrac{d}{d\epsilon}\Big|_{\epsilon=0}(\varphi(t;v_{0}+\epsilon\delta
v_{0}),v(t;v_{0}+\epsilon\delta v_{0}))
\]
which is a continuous curve in
$H^\ell_\C\times H^\ell_\C$. Differentiating \eqref{e:crucial} in
direction $\delta v_0$ in $H^\ell_\C$ at $v_0$ in $V_{\ell,\C}^{(k)}$
yields
\begin{eqnarray}
(A_k-1)(\delta\varphi(t))-\varphi'(t)\left(\int_0^t(\varphi'(\tau))^{2k-3}\,d\tau\right) A_k(\delta v_0)=\nonumber\\
=P_{2k-1}(\varphi(t),v(t);\delta\varphi(t),\delta v(t))\label{e:crucial1}
\end{eqnarray}
where $P_{2k-1} :U_{\ell,\C}^{(k)}\times H^\ell_\C\times H^\ell_\C\times H^\ell_\C\to H^{\ell-2k+1}_\C$ is analytic.

\begin{Prop}\hspace{-2mm}{\bf .}\label{Proposition 7.1} 
Let $k \geq 1$, $\ell_k=2k+2$, and $\ell=\ell_{k}+n$, $n\ge0$. 
Assume that $v_0$ is in $V^{(k)}_{\ell,\C}$. Then
\[ 
(d_{v_0}\mathrm{Exp}^\C_{k,\ell})(H^\ell_\C\backslash H^{\ell+1}_\C)
\subseteq H^\ell_\C\backslash H^{\ell+1}_\C\,.
\]
\end{Prop}
{\it Proof.} Assume $v_0\in V^{(k)}_{\ell,\C}$, and let
$(\varphi(\cdot;v_0),v(\cdot;v_0))$ in
$C^{1}((-2,2),U_\ell^{(k)}\times H^\ell_\C)$ be the unique solution of
\eqref{ds-complex} issuing from $(\mathrm{id},v_{0})$ as guaranteed by
Theorem~\ref{Th:extension}. As before, it follows from
\eqref{U_1-complex} that the factor in front of $A_{k}\delta v_0$ in
formula \eqref{e:crucial1}, evaluated at $t=1$, is a non-zero function
in $H^{\ell-1}_{\C}\subseteq H^{\ell-2k+1}_\C$ whereas the term on the r.h.s. of this identity is in $H^{\ell-2k+1}_\C$. Hence, just as in the proof
of Proposition~\ref{Proposition 6.1} which followed from analyzing
\eqref{e:crucial}, the statement of Proposition~\ref{Proposition 7.1}
can be obtained from \eqref{e:crucial1}, evaluated at $t=1$.
\finishproof

\appendix

\section{Analytic maps between Fr{\'e}chet
  spaces}\label{Appendix:Frechet_spaces}
For the convenience of the reader we collect in this appendix some
definitions and notions from the calculus in Fr{\'e}chet spaces and
present an inverse function theorem valid in a set-up for Fr\'echet
spaces which is suitable for our purposes. For more details on the
theory of smooth functions in Fr{\'e}chet spaces we refer the reader
to \cite{Ham}. For the theory of analytic functions in Fr{\'e}chet
spaces, we follow the approach developped in \cite{BS1,BS2} (cf. also
\cite{KMich}). In the sequel $\K$ denotes either the field $\C$ of complex
numbers or the field $\R$ of real numbers.

\medskip

{\em Fr{\'e}chet spaces:} Consider the pair $(X,\{||\cdot||_n\}_{n\ge0})$ where $X$ is 
a vector space over $\K$ and $\{||\cdot||_n\}_{n\ge0}$ is a countable collection of seminorms.  A
topology on $X$ is defined in the usual way as follows: A basis of
open neighborhoods of $0$ in $X$ is given by the sets
\[
U_{\epsilon,k_1,...,k_s}\eqdef\{x\in X\;:\;||x||_{k_j}<\epsilon\;\;\;\forall 1\le j\le s\}
\]
where $s,k_1,...,k_s$ are nonnegative integers and $\epsilon>0$.  Then
the topology on $X$ is defined as the collection of open sets
generated by the sets $x+U_{\epsilon,k_1,...,k_s}$, for arbitrary $x$
in $X$ and arbitrary $s,k_1,...,k_s$ in $\Z_{\ge0}$ and $\epsilon>0$.
In this way, $X$ becomes a topological vector space. Note that a
sequence $(x_k)_{k\ge0}$ converges to $x$ in $X$ iff, for any $n\ge
0$, $||x_k-x||_n\to 0$ as $k\to +\infty$.

\medskip

Moreover, the topological vector space $X$ described above is {\em
  Hausdorff} iff, for any $x$ in $X$, 
$||x||_n=0$ for every $n$ in $\Z_{\ge0}$ implies $x=0$. A sequence $(x_k)_{k\in\N}$ is called
{\em Cauchy} iff it is a Cauchy sequence with respect to any of the
seminorms. By definition, $X$ is complete iff every Cauchy sequence
converges in $X$.

\begin{Def}\hspace{-2mm}{\bf .}
A pair $(X,\{||\cdot||_n\}_{n\ge0})$ consisting of a topological
vector space $X$ and a countable system of seminorms
$\{||\cdot||_n\}_{n\ge0}$ is called a {\em Fr{\'e}chet space}
\footnote{Unlike for the standard notion of a Fr{\'e}chet space,
here the countable system of seminorms defining the topology of
$X$ is part of the structure of the space.} iff the topology of
$X$ is the one induced by $\{||\cdot ||_n\}_{n \in\Z_\ge0}$, and $X$
is Hausdorff and complete.
\end{Def}

The space of continuous maps $f:U\to Y$ from an open subset
$U\subseteq X$ into the Fr\'echet space $Y$ is denoted by $C^{0}(U,Y)$.

\medskip

{\em $C_{F}^{1}$-differentiability:}
Let $f : U\subseteq X\to Y$ be a map from an open set $U$ of a Fr{\'e}chet space $X$
to a Fr{\'e}chet space $Y$.

\begin{Def}\hspace{-2mm}{\bf .}\label{Def:directiona_derivative}
If the limit 
\[
  \lim_{\epsilon\in\K,\epsilon\to 0}\frac{1}{\epsilon}(f(x+\epsilon h)-f(x))
\]
in $Y$ exists with respect to the Fr\'echet topology of $Y$, we say that {\em
  $f$ is differentiable at $x$ in the direction $h$}. The limit is
declared to be the {\em directional derivative} of $f$ at the point
$x$ in $U$ in the direction $h$ in $X$. Following \cite{BS1,BS2}, we
denote it by $\delta_xf(h)$.
\end{Def}

\begin{Def}\hspace{-2mm}{\bf .}\label{Def:C^{1}}
If the directional derivative $\delta_xf(h)$ exists for any $x$ in $U$ and
any $h$ in $X$, and the map
\[
(x,h)\mapsto \delta_xf(h),\;U\times X\to Y
\]
is continuous with respect to the Fr{\'e}chet topology on $U\times X$
and $Y$, then $f$ is called {\em continuously differentiable on $U$}.
The space of all such maps is denoted by
$C_{F}^{1}(U,Y)$.\footnote{Note that even in the case where $X$ and
$Y$ are Banach spaces this definition of continuous
differentiability is weaker than the usual one (cf. \cite{Ham}). In
order to distinguish it from the classical one we write $C_{F}^{1}$
instead of $C^{1}$. We refer to \cite{Ham} for a discussion of the
reasons to introduce the notion of $C_{F}^{1}$-differentiability.}
A map $f : U\to V$ from an open set $U\subseteq X$ onto an open set
$V\subseteq Y$ is called {\em a $C_{F}^{1}$-diffeomorphism} if $f$ is
a homeomorphism and $f$ as well as $f^{-1}$ are $C_{F}^{1}$-smooth.
\end{Def}

\begin{Lemma}\hspace{-2mm}{\bf .}\label{Rem:no_kernel}
Let $U\subseteq X$ be an open subset. Then 
\begin{itemize}
\item[(i)] $C^{1}_{F}(U,Y)\subseteq C^{0}(U,Y)$
\item[(ii)]Assume that a map $f$ in $C^{1}(U,Y)$ is a
  $C^{1}_{F}$-diffeomorphism onto an open subset $V\subseteq Y$. Then,
  for any $x$ in $U$, $ \delta_xf : X\to Y $ is a linear isomorphism.
\end{itemize}  
\end{Lemma} {\em Proof.} Statement $(i)$ follows from \cite[Theorem
3.2.2.]{Ham} whereas statement $(ii)$ is a consequence of
\cite[Theorem 3.3.4.]{Ham}.\finishproof

\medskip

{\em Analytic functions in Fr{\'e}chet spaces:} Let $X$ and $Y$ be
Fr{\'e}chet spaces over $\K$ and let $f : U\subseteq X\to Y$ be a map
from an open set $U\subseteq X$ into $Y$. A map $f_s : X\to Y$ is
called a {\em homogeneous polynomial} of degree $s\in\Z_{\ge 0}$ if
there exists a $s$-linear symmetric map $\mathtt{f}_s : X^s\to Y$ such
that $f_s(x)=\mathtt{f}_s(x,...,x)$ for any $x$ in $X$ (cf.
\cite[Definition 2]{BS1}).

\begin{Def}
\label{Def:analytic_functions}
Following \cite[Definition 5.6]{BS2}, a continuous function $f:U\to Y$
is called {\em analytic} if, for any $x$ in $U$, there exist an open
neighborhood $V$ of $0$ in $X$ and a sequence of continuous
homogeneous polynomials $(f_s)_{s\ge 0}$, $\deg f_s=s$, such that
$x+V\subseteq U$ and, for any $h$ in $V$, $ f(x+h)=\sum_{s=0}^\infty
f_s(h)$ converges in $Y$\footnote{In case $\K=\R$, an analytic
  function $f:U\to Y$ is sometimes called real analytic.}.
\end{Def}
We will need the following lemma.
\begin{Lemma}\hspace{-2mm}{\bf .}\label{Lem:useful}
Assume that $X$ and $Y$ are Fr{\'e}chet spaces over $\C$, $U\subseteq
X$ is an open subset of $X$, and $f:U\to Y$ is in $C^{1}_{F}(U,Y)$. Then $f$ is analytic.
\end{Lemma} 

\begin{Rem}
\label{RemarkA.6bis}
The converse of Lemma~\ref{Lem:useful} is true as well. More
precisely, assume that $X$ and $Y$ are $\K$-Fr\'echet spaces and
$f:U\to Y$ is analytic. Then, by the definition of an analytic map,
$f$ is in $C^{n}_{F}(U,Y)$ for any $n\ge0$.
\end{Rem}

{\em Proof of Lemma~\ref{Lem:useful}.} According to \cite[Theorem
6.2]{BS2}, it suffices to prove that $f$ is continuous and that it is
analytic on affine lines. By Lemma~\ref{Rem:no_kernel} $(i)$, $f$ is
continuous. By Definition~\ref{Def:C^{1}}, it follows that, for any $x$
in $U$, and any $h$ in $X$, the map $ f_{x,h} : z\mapsto f(x+zh)$ with
values in $Y$, defined on the open set $\{z\in\C\;:\;x+zh\in
U\}\subseteq\C$ is (complex) differentiable. In particular, by
\cite[Theorem 3.1]{BS2}, $f_{x,h}(z)$ is analytic.  Hence, $h$ in $X$
being arbitrary, $f$ is analytic on affine lines.\finishproof

\medskip

{\em Analytic functions in Fr{\'e}chet spaces over $\R$:} Now, assume
that $X$ and $Y$ are Fr{\'e}chet spaces over $\R$. Denote by
$X_\C=X\otimes\C$ the complexification of $X$.  The following theorem
follows directly from \cite[Theorem 3]{BS1} and \cite[Theorem
7.1]{BS2}.
\begin{Th}\hspace{-2mm}{\bf .}
  Let $U$ be an open subset of $X$. A function $f : U\to Y$
  is analytic iff there exists a complex neighborhood
  ${\tilde U}\supseteq U$ in $X_\C$ and an analytic function ${\tilde
    f} : {\tilde U}\to Y_\C$ such that ${\tilde f}|_U=f$.
\end{Th}
In this paper we consider mainly the following spaces:

\medskip

{\em Fr{\'e}chet space $C^\infty(\T)$:}
The space $C^\infty(\T)\equiv C^\infty(\T,\R)$ denotes the real vector space
of real-valued $C^\infty$-smooth, $1$-periodic functions $u :\R\to\R$.
The topology on $C^\infty(\T)$ is induced by the countable system of Sobolev norms:
\[
\|
u\|_n\eqdef\|u\|_{H^{n}}=\Big(\sum_{j=0}^n\int_0^1[\p_{x}^{j}u(x)]^2\;dx\Big)^{1/2},\quad
n\ge 0.
\]

\medskip

{\em Fr{\'e}chet manifold $\D$:} 
By definition, $\D$ denotes the group of $C^\infty$-smooth positively oriented
diffeomorphisms of the torus $\T=\R/\Z$.  A Fr{\'e}chet manifold
structure on $\D$ can be introduced as follows: Passing to the universal cover $\R\to\T$,
any element $\varphi$ of $\D$ gives rise to a smooth diffeomorphism of $\R$ in $C^\infty(\R,\R)$,
again denoted by $\varphi$, satisfying the {\em normalization condition}
\[
-1/2<\varphi(0)<1/2\quad\quad\mbox{or}\quad\quad 0<\varphi(0)<1.
\]
The function $f\eqdef\varphi-\mathrm{id}$ is $1$-periodic and hence
lies in $C^\infty (\T)$. Moreover $f'(x)>-1$ for any $x$ in $\R$.  The
above normalizations give rise, respectively, to two charts
$\mathfrak{U}_1$, $\mathfrak{U}_2$ of $\D$ with $\mathfrak{U}_1\cup
\mathfrak{U}_2=\D$, defined by
\[
\mathfrak{F}_j : \mathfrak{V}_j\to \mathfrak{U}_j, \quad
f\mapsto\varphi:=\mathrm{id}+f
\]
where $j=1,2$, and 
\begin{eqnarray}
\mathfrak{V}_1&\eqdef&\{f\in C^\infty(\T)\;|\;-1/2<f(0)<1/2\;\mbox{and}\;f'>-1\}\nonumber\\
\mathfrak{V}_2&\eqdef&\{f\in C^\infty(\T)\;|\;0<f(0)<1\;\mbox{and}\;f'>-1\}.\nonumber
\end{eqnarray}
As $\mathfrak{V}_1, \mathfrak{V}_2$ are both open subsets in the Fr{\'e}chet space $C^\infty(\T)$, the
construction above gives an atlas of Fr{\'e}chet charts of $\D$.
In this way, $\D$ is a Fr{\'e}chet manifold modeled on $C^\infty(\T)$.

\medskip

{\em Hilbert manifold $\D^s(\T)$ ($s\ge 2$):}
$\D^{s}=\D^s(\T)$ denotes the group of positively oriented
bijective transformations of $\T$ of class $H^{s}$.  By definition, a
bijective transformation $\varphi$ of $\T$ is of class $H^{s}$ iff the
lift ${\tilde\varphi} : \R\to\R$ of $\varphi$, determined by the
normalization, $0\le{\tilde\varphi}(0)<1$, and its inverse
${\tilde\varphi}^{-1}$ both lie in the Sobolev space
$H^s_{loc}(\R,\R)$.  As for $\D$ one can introduce an atlas for
$\D^s$ with two charts in $H^s$, making $\D^s$ a Hilbert manifold
modeled on $H^s$.

\medskip

{\em Hilbert approximations:}
Assume that for a given Fr{\'e}chet space $X$ over $\K$ there is a sequence of $\K$-Hilbert spaces
$(X_n,||\cdot||_n)_{n\ge0}$ such that
\[
X_0\supseteq X_1\supseteq X_2\supseteq...\supseteq
X\;\;\;\mbox{and}\;\;\;X=\displaystyle\bigcap_{n =0}^\infty X_n
\]
where $\{||\cdot||_n\}_{n\ge0}$ is a sequence of norms inducing the
topology on $X$ so that $||x||_0\le||x||_1\le||x||_2\le...$ for any
$x$ in $X$.  Such a sequence of Hilbert spaces $(X_n,||\cdot||_n)_{n\ge0}$ is called 
a {\em Hilbert approximation} of the Fr{\'e}chet space $X$.  
For Fr{\'e}chet spaces admitting Hilbert approximations one can prove 
the following version of the inverse function theorem.

\begin{Th}\hspace{-2mm}{\bf .}\label{Th:IFT} 
Let $X$ and $Y$ be Fr\'echet spaces over $\K=\C$ or $\R$ with Hilbert approximations 
$(X_n, \|\cdot\|_n)_{n\geq 0}$, and respectively $(Y_n,|\cdot|_n)_{n\geq 0}$. 
Let $f : V_0 \rightarrow U_0$ be a map between the open subsets $V_0 \subseteq X_0$ and $U_0 \subseteq Y_0$ 
of the Hilbert spaces $X_0$, respectively $Y_0$. Define, for any $n \geq 0$,
\[
V_n:= V_0 \cap X_n ,\quad\quad\quad\quad U_n := U_0 \cap Y_n .
\]
Assume that, for any $n\ge0$, the following properties are satisfied:
\begin{itemize}
\item[{\rm (a)}] $f : V_0 \rightarrow U_0$ is a bijective $C^{1}$-map,
and, for any $x$ in $V:=V_0\cap X$, $d_x f : X_0 \rightarrow Y_0$ is a linear isomorphism;
\item[{\rm (b)}] $f(V_n) \subseteq Y_n$, and the restriction 
$f\big\arrowvert_{V_n} : V_n \rightarrow Y_n$ is a $C^{1}$-map;
\item[{\rm (c)}] $f(V_n)\supseteq U_n$;
\item[{\rm (d)}] for any $x$ in $V$, $d_x f(X_n \backslash X_{n + 1})\subseteq Y_n \backslash Y_{n + 1}$. 
\end{itemize}   
Then for the open subsets $V:=V_{0}\cap X\subseteq X$ and $U:= U_0 \cap Y \subseteq Y$, one has $f(V)\subseteq U$ and the map 
$f_\infty := f\big\arrowvert_V : V \rightarrow U$ is a $C^{1}_{F}$-diffeomorphism.
\end{Th} {\it Proof.} By properties (a) and (b), $f_n:=
f\big\arrowvert _{V_n} : V_n \rightarrow U_n$ is a well-defined,
injective $C^{1}$-map. By (c), $f_n$ is onto. Hence, $f_\infty :=
f\big\arrowvert_V : V \rightarrow U$ is bijective. In order to prove
that $f_\infty : V\rightarrow U$ is a $C^{1}_{F}$-diffeomorphism,
consider, for any $n\ge0$, and any $x$ in $V$, the differential $d_x f_n
: X_n \rightarrow Y_n$.  As $(d_x f)\big\arrowvert _{X_n} = d_x f_n$
and, by (a), $d_x f : X_0 \rightarrow Y_0$ is bijective, one concludes
that $d_x f_n$ is one-to-one. We prove by induction (with respect to
$n$) that, for any $x$ in $V$, $d_x f_n : X_n \rightarrow Y_n$ is onto.
For $n = 0$ ($V\subseteq V_{0}$), the statement is true by property
(a).  Next, assume that for arbitrary positive integer $n$, and
arbitrary $x$ in $V$, $d_x f_{n - 1} : X_{n - 1} \rightarrow Y_{n -
  1}$ is onto.  Then, for any $x$ in $V$, and $\eta$ in $Y_{n}\subseteq
Y_{n-1}$, there exists a (unique) $\xi$ in $X_{n-1}$ verifying $d_x
f_{n - 1}(\xi ) = \eta $. By property (d), it follows that $\xi$
belongs to $X_{n}$.  In other words, for any given $n\ge0$, and any $x$
in $V$, we have that the map $d_x f_n : X_n \rightarrow Y_n$ is
bijective, and thus, by Banach's theorem, the inverse $(d_x f_n)^{-1}:
Y_n \rightarrow X_n$ is a bounded linear operator. As, for any $n\ge0$,
$f_{n}$ is $C^{1}$-smooth, the map
\begin{equation}\label{A.20} 
V'_n\times X_n\rightarrow Y_n,\,(x,\xi)\mapsto d_x f_n(\xi )
\end{equation}
is continuous and, by the inverse function theorem it follows that
\begin{equation}\label{A.21} 
U'_n\times Y_n\rightarrow X_n,\,(y,\eta)\mapsto d_y (f^{-1}_n)(\eta )
\end{equation}
is continuous as well. Here $V'_n$ ($U'_n$) denotes the subset $V$
($U$) with the topology induced by $|\cdot|_n$ ($\|\cdot\|_n$). As for
any $x$ in $V$, and $n\geq 0$,
\[ 
\delta_x f_\infty = d_x f_n \big\arrowvert_X
\]
one gets from \eqref{A.20} - \eqref{A.21} that
\[ 
V \times X \rightarrow Y , \ (x, \xi ) \rightarrow \delta_x f_\infty(\xi )
\]
and
\[ 
U \times Y \rightarrow X , \ (x, \eta ) \mapsto \delta_y f^{-1}_\infty (\eta )
\]
are continuous. In particular, one concludes (cf. Definition~\ref{Def:C^{1}}) that
\[ 
f_\infty : V \rightarrow U
\]
is a $C^{1}_{F}$-diffeomorphism.  
\finishproof

\bibliographystyle{plain}

\end{document}